\documentclass{gtart}

\def\ifplaintex{\expandafter\ifx\csname documentclass\endcsname\relax}

\ifplaintex 
\else
\fi

\expandafter\ifx\csname beginpicture\endcsname\relax
\expandafter\ifx\csname documentclass\endcsname\relax
\input pictex \else
\input prepictex \input pictex \input postpictex \fi\fi

\def\gt{{\mathsurround=0pt\it $\cal G\mskip-2mu$eometry \&\ 
$\cal T\!\!$opology}}        

\def\gtp{{\mathsurround=0pt\it $\cal G\mskip-2mu$eometry \&\ 
$\cal T\!\!$opology $\cal P\!$ublications}}  

\def\lognumber#1{\def\thelognumber{#1}}
\def\volumenumber#1{\def\thevolumenumber{#1}}
\def\papernumber#1{\def\thepapernumber{#1}}
\def\volumeyear#1{\def\thevolumeyear{#1}}

\def\pagenumbers#1#2{\def\startpage{#1}\def\finishpage{#2}}
\def\published#1{\def\publishdate{#1}}
\def\proposed#1{\def\theproposer{#1}}
\def\seconded#1{\def\theseconders{#1}}
\def\received#1{\def\receiveddate{#1}}
\def\revised#1{\def\reviseddate{#1}}
\def\accepted#1{\def\accepteddate{#1}}

\let\\\par\let\thelognumber\relax
\let\thevolumenumber\relax\let\thepapernumber\relax
\let\thevolumeyear\relax\let\thesamplenumber\relax\let\startpage\relax
\let\finishpage\relax\let\publishdate\relax\let\receiveddate\relax
\let\reviseddate\relax\let\accepteddate\relax\let\theasciititle\relax
\let\theasciiauthors\relax
\let\theasciiabstract\relax
\let\theasciiemail\relax\let\theshortauthors\relax\let\theshorttitle\relax

\long\def\maketitlep{   

\count0=\startpage

\gt\hfill      
\beginpicture
\setcoordinatesystem units <0.33truein, 0.33truein> point at 2.2 0.9
\setplotsymbol ({$\cal G$})
\plotsymbolspacing=9truept
\circulararc 315 degrees from 0 1 center at 0 0
\setplotsymbol ({$\cal T$})
\circulararc 315 degrees from 1 -1 center at 1 0
\endpicture
\break
{\small\ifx\thesamplenumber\relax 
Volume \else Sample
\fi\thevolumenumber\ (\thevolumeyear)
\startpage--\finishpage\nl
Published: \publishdate}
\vglue 0.5truein plus 0.4fil minus 0.1truein

{\parskip=0pt\leftskip 0pt plus 1fil\def\\{\par\smallskip}{\ifplaintex\large
\else\Large\fi\bf\thetitle}\par\medskip}   

\vglue 0pt plus 0.1fil 

{\parskip=0pt\leftskip 0pt plus 1fil\def\\{\par}{\sc\theauthors}
\par\medskip}

\vglue 0pt plus 0.1fil 

{\small\parskip=0pt\let\newline\\
{\leftskip 0pt plus 1fil\def\\{\par}{\sl\theaddress}\par}
\expandafter\ifx\theemail\relax    
\relax\else\vglue 5pt plus 0.02fil minus 2pt\def\\{\stdspace{\rm 
and}\stdspace} 
\cl{Email:\stdspace\tt\theemail}\fi
\ifx\theurl\relax                  
\relax\else\vglue 5pt plus 0.02fil minus 2pt\def\\{\stdspace{\rm 
and}\stdspace}
\cl{URL:\stdspace\tt\theurl}\fi\par}

\vglue 7pt plus 0.3fil minus 3pt

{\bf Abstract}
\vglue 5pt plus 0.1fil minus 2pt

\theabstract

\vglue 7pt plus 0.3fil minus 3pt

{\bf AMS Classification numbers}\quad Primary:\quad \theprimaryclass

Secondary:\quad \thesecondaryclass

\vglue 5pt plus 0.3fil minus 2pt

{\bf Keywords}\quad \thekeywords

\vglue 10pt plus 0.5fil minus 5pt

{\small  Proposed: \theproposer\hfill Received: \receiveddate\nl
Seconded: \theseconders\hfill 
\ifx\reviseddate\relax                         
Accepted: \accepteddate                        
\else
Revised: \reviseddate                          
\fi}
\eject
}       

\let\maketitlepage\maketitlep
\let\maketitle\maketitlepage

\font\phead=cmsl9 scaled 950
\font=cmsl9 scaled 1050
\font\pnum=cmbx10 scaled 913
\font\lnum=cmbx10 
\font\pfoot=cmsl9 scaled 950
\font=cmsl9 scaled 1050
\ifplaintex
\headline{\vbox to 0pt{\vskip -4.5mm\line{\small\phead\ifnum
\count0=\startpage ISSN 1364-0380 (on line)
1465-3060 (printed) \hfill {\pnum\folio}\else\ifodd\count0\def\\{ }%
\ifx\theshorttitle\relax\thetitle\else\theshorttitle\fi\hfill{\pnum\folio}
\else\def\\{ and }{\pnum\folio}\hfill\ifx\theshortauthors\relax\theauthors
\else\theshortauthors\fi\fi\fi}\vss}}
\footline{\vbox to 0pt{\vglue 0mm\line{\small\pfoot\ifnum\count0=\startpage
\copyright\ \gtp\hfill\else
\gt, Volume \thevolumenumber\ (\thevolumeyear)\hfill\fi}\vss
}}
\else
\makeatletter
\def\@oddhead{{\small\lhead\ifnum\count0=\startpage ISSN 1364-0380 (on line)
1465-3060 (printed) \hfill {\lnum\number\count0}\else\ifodd\count0
\def\\{ }\ifx\theshorttitle\relax \thetitle \else\theshorttitle\fi\hfill
{\lnum\number\count0}\else\def\\{ and }{\lnum\number\count0}
\hfill\ifx\theshortauthors\relax 
\theauthors\else\theshortauthors\fi\fi\fi}}\def\@evenhead{\@oddhead}
\def\@oddfoot{\small\lfoot\ifnum\count0=\startpage\copyright\ \gtp\hfill\else
\gt, Volume \thevolumenumber\ (\thevolumeyear)\hfill\fi}
\def\@evenfoot{\@oddfoot}
\makeatother
\fi

\newwrite\gtoutfile
\long\gdef\makeheadfile{  
{\def\\{, }\def\s{ }
\immediate\openout\gtoutfile head.xxx
\immediate\write\gtoutfile{Proxy-for: \ifx\theasciiauthors\relax
\theauthors\else\theasciiauthors\fi\s<\ifx\theasciiemail\relax\theemail\else\theasciiemail\fi>}
\immediate\write\gtoutfile{\noexpand\\}
\immediate\write\gtoutfile{Authors: \ifx\theasciiauthors\relax
\theauthors\else\theasciiauthors\fi}
{\def\\{ }\immediate\write\gtoutfile{Title: \ifx\theasciititle\relax
\thetitle\else\theasciititle\fi}}
\immediate\write\gtoutfile{Subj-class: GT or SG or MG etc}
\immediate\write\gtoutfile{MSC-class: \theprimaryclass\ifx\thesecondaryclass\relax\else, \thesecondaryclass\fi}
\immediate\write\gtoutfile{Journal-ref: Geom. Topol. \thevolumenumber
(\thevolumeyear) \startpage-\finishpage}
\immediate\write\gtoutfile{Comments: Published by Geometry and Topology at}
\immediate\write\gtoutfile{\s\s http://www.maths.warwick.ac.uk/gt/GTVol\thevolumenumber/paper\thepapernumber.abs.html}
\immediate\write\gtoutfile{\noexpand\\}
\immediate\write\gtoutfile{}
\ifx\theasciiabstract\relax
\immediate\write\gtoutfile{\theabstract}\else
\immediate\write\gtoutfile{\theasciiabstract}\fi
\immediate\write\gtoutfile{}
\immediate\write\gtoutfile{\noexpand\\}
\immediate\write\gtoutfile{}
\immediate\closeout\gtoutfile}}  

\def\maketitlepage{\maketitlep\makeheadfile}
\let\maketitle\maketitlepage

\lognumber{258}
\volumenumber{8}\papernumber{3}\volumeyear{2004}
\pagenumbers{77}{113}
\received{29 May 2002}
\revised{21 January 2004}
\published{22 January 2004}
\accepted{13 December 2003}
\proposed{Cameron Gordon}
\seconded{David Gabai, Joan Birman}

\usepackage{epsfig}
\usepackage{psfrag}
\usepackage{amssymb}
\usepackage{amsmath}
\let\url\undefined
\usepackage{hyperref}

\newcommand{\ie}{{i.e.}}
\newcommand{\etal}{{et.\thinspace al.}}

\newcommand{\cf}{{\mathcal{F}}}

\renewcommand{\setminus}{{\smallsetminus}}

\newcommand{\RR}{{\mathbb{R}}}

\newcommand{\DD}{{\mathbb{D}}}

\newcommand{\PP}{{\mathbb{P}}}

\newcommand{\from}{{\colon}} 

\newcommand{\double}[1]{{\widetilde{#1}}}

\newcommand{\closure}[1]{{\overline{#1}}}

\newcommand{\euler}{{\chi}} 
\newcommand{\homeo}{{\medspace \cong \medspace}} 
 
\newcommand{\neigh}{{\eta}} 
\newcommand{\frontier}{{\operatorname{fr}}} 
\newcommand{\fr}{{\operatorname{fr}}} 
\newcommand{\bdy}{{\partial}} 
\newcommand{\disjoint}{{\medspace \amalg \medspace}}
\newcommand{\rel}{{\operatorname{rel}}}
 
\newcommand{\cross}{{\times}}

\newcommand{\interior}{{\operatorname{interior}}}

\theoremstyle{plain}
\newtheorem{theorem}{Theorem}[section]

\newtheorem{lemma}[theorem]{Lemma}
\newtheorem{conjecture}[theorem]{Conjecture}
\newtheorem{claim}[theorem]{Claim}
\newtheorem*{proofclaim}{Claim}
\newtheorem{proposition}[theorem]{Proposition}

\theoremstyle{definition}
\newtheorem*{define}{Definition}
\newtheorem{remark}[theorem]{Remark}

\newtheorem*{fact}{Fact}
\newtheorem*{question}{Question}

\newsavebox{\savepar}

\newcommand{\blockconst}{{32}}
\newcommand{\dcpconst}{{\exp(2^{16} t^2)}}

\newcommand{\blocky}[1]{{\widehat{#1}}}
\newcommand{\Mobius}{{M\"obius}}

\begin{document}

\title{The disjoint curve property}
\author{Saul Schleimer}
\address{Department of Mathematics, UIC\\851 South Morgan 
Street\\Chicago, Illinois 60607, USA}
\email{\href{mailto:saul@math.uic.edu}{saul@math.uic.edu}}
\urladdr{\href{http://www.math.uic.edu/~saul}{http://www.math.uic.edu/\char'176saul}}

\begin{abstract}
A Heegaard splitting of a closed, orientable three-manifold satisfies
the disjoint curve property if the splitting surface contains an
essential simple closed curve and each handlebody contains an
essential disk disjoint from this curve [Thompson, 1999].  A splitting
is full if it does not have the disjoint curve property.  This paper
shows that in a closed, orientable three-manifold all splittings of
sufficiently large genus have the disjoint curve property.  From this
and a solution to the generalized Waldhausen conjecture it would
follow that any closed, orientable three manifold contains only
finitely many full splittings.
\end{abstract}
\maketitlepage

\section{History and overview}
\label{Sec:Intro}


The classification of Heegaard splittings of three-manifolds is a long
outstanding problem.  It is natural to add restrictions to the class
of three-manifolds considered.  Haken's lemma~\cite{Haken68}, that all
splittings of a reducible manifold are themselves reducible, can be
considered one of the first results in this direction.  {\em Strong
  irreducibility} was introduced by Casson and
Gordon~\cite{CassonGordon87} as a generalization of irreducibility.
They showed that all splittings of a non-Haken manifold are either
reducible or strongly irreducible.  Thompson~\cite{Thompson99} later
defined the {\em disjoint curve property} (DCP) as a further
generalization of reducibility.  She deduced that all splittings of a
toroidal three-manifold have the disjoint curve property.
Hempel~\cite{Hempel01}, using the classification of splittings of
Seifert fibred spaces~\cite{MoriahSchultens98}, has shown that each of
these has the disjoint curve property.  Thus, in any three-manifold
which is reducible, toroidal, or a Seifert fibred space, all Heegaard
splittings have the disjoint curve property.

However, it is certainly not the case that all splittings of all
manifolds have the disjoint curve property.  Hempel (referring to Luo)
in~\cite{Hempel01} adapts an argument of Kobayashi~\cite{Kobayashi88b}
to produce examples of splittings which are {\em full} and are in fact
arbitrarily far from having the DCP.

The remainder of this section outlines the paper and the results
contained herein.  The required ideas from the theories of Heegaard
splittings and normal surfaces are laid out in
Sections~\ref{Sec:Splittings} and~\ref{Sec:Normal}.  Next, {\em
  blocks} and various associated objects are defined in
Section~\ref{Sec:Blocks}.  These are used in the proof of the main
technical proposition:

\medskip
{\bf Proposition~\ref{Prop:AlmostNormalPlusHighGenusImpliesDCP}}\qua\sl
If $H \subset (M,T_M)$ is an almost normal Heegaard splitting with
sufficiently large genus then $H$ has the disjoint curve property.
\rm\medskip

The rough idea is that a high genus almost normal surface must
decompose as a {\em Haken sum}.  The {\em exchange annuli} of this sum
should then give the desired pair of disks and disjoint curve.  This
rough idea is not quite right, of course, and the proof of the
proposition occupies the bulk of the paper.  An outline is given
in Section~\ref{Sec:MainProposition}.  

Now, not every Heegaard splitting is given as an almost normal
surface.  However, it is a result of Rubinstein and
Stocking~\cite{Stocking00} that every strongly irreducible splitting
may be normalized.  (Stated as Theorem~\ref{Thm:NormalizeStrIrred}
below.)  Thus our main result, Theorem~\ref{Thm:HighGenusImpliesDCP},
is a corollary of
Proposition~\ref{Prop:AlmostNormalPlusHighGenusImpliesDCP}:

\medskip
{\bf Theorem~\ref{Thm:HighGenusImpliesDCP}}\qua  
{\sl 
Every splitting in $M$ of sufficiently large genus has the DCP.
}
\medskip

That this result is sharp in general is shown by the Casson-Gordon-Parris
examples:

\begin{fact} 
There is a closed, Haken three-manifold containing strongly
irreducible splittings of arbitrarily large genus.
(See~\cite{Kobayashi92} or~\cite{Sedgwick97}.)
\end{fact}

Finally, it would be remiss not to mention that Jaco and
Rubinstein~\cite{JacoRubinstein02} have claimed a solution to the
generalized Waldhausen conjecture:

\begin{conjecture}
\label{Conj:GeneralizedWaldhausen}
If $M$ is closed, orientable, and atoroidal then $M$ contains only
finitely many strongly irreducible splittings, up to isotopy, in each
genus.
\end{conjecture}

This, when combined with Theorem~\ref{Thm:HighGenusImpliesDCP}, would
give:

\medskip
{\bf Conjecture~\ref{Conj:OnlyFinManyFullSplittings}}\qua
{\sl 
In any closed, orientable three-manifold there are only finitely many
full Heegaard splittings, up to isotopy.
}
\medskip

The paper ends in Section~\ref{Sec:Conjecture} by explaining this and
other, related, open questions.

I thank Kevin Hartshorn and David Bachman for many important
conversations in which the main question of the paper was formulated.
I also thank Yoav Moriah for his valuable comments on an early version
of this paper.  I owe a further debt of gratitude to my thesis
advisor, Andrew Casson.

\section{Heegaard splittings}
\label{Sec:Splittings}

This section recalls standard notation and reviews several notions
from the theory of Heegaard splittings.  (For an excellent survey of
the subject, see~\cite{Scharlemann02}.)  The slightly nonstandard
notions of {\em disjoint} and {\em joined} essential surfaces are
given as well as the concept of a {\em full splitting}.

In this paper, $M$ denotes a closed, orientable three-manifold.
Typically $N$ will denote a compact, orientable three-manifold with
possibly non-empty boundary.  Recall that a three-manifold is {\em
reducible} if it contains an embedded two-sphere which does not bound
an embedded three-ball in $M$.

A {\em handlebody} is a compact three-manifold which is homeomorphic
to a closed regular neighborhood of a finite, polygonal, connected
graph embedded in $\RR^3$.  A {\em Heegaard splitting} (or simply a
{\em splitting}), $H \subset M$, is a embedded, orientable,
separating, closed surface such that the closure of each component of
$M \setminus H$ is a handlebody.  These are typically denoted $V$ and
$W$.  The {\em genus} of the splitting is the genus of $H$, $g(H) = 1
-\euler(V)$.

Suppose $A \subset V$ is a connected and properly embedded surface.
Let $\double{A}$ be the frontier of a closed regular neighborhood of
$A$, taken in $V$. Then the surface $\double{A}$ is the {\em double}
of $A$.

\begin{define}
If $D$ is a disk properly embedded in a handlebody $V$ then $D$ is
{\em essential} if $\bdy D$ is not null-homotopic in $\bdy V$.  If $A$
is an annulus properly embedded in $V$ then $A$ is {\em essential} if
$A$ is incompressible in $V$ and not boundary-parallel ($\rel~\bdy A$)
to an annulus embedded in $\bdy V$.  A properly embedded \Mobius~band
$A \subset V$ is {\em essential} exactly when $\double{A}$ is.
(Recall that a surface $F$ properly embedded in a compact manifold $N$
is {\em incompressible} is for every embedded disk $(D, \bdy D)
\subset (N, F)$ with $D \cap F = \bdy F$ there is an embedded disk
$(E, \bdy E) \subset (F, \bdy D)$.  Also, $F \subset N$ is {\em
  boundary-parallel} if there is a subsurface $G \subset \bdy N$ with
$\bdy F = \bdy G$, $F \cup G$ bounds a submanifold of $N$ homeomorphic
to the product $F \cross I$, and where $F \cross \{1\} = F$, $(F
\cross \{0\}) \cup (\bdy F \cross I) = G$.)
\end{define}

\begin{remark}
\label{Rem:Essential}
This nonstandard definition of ``essential'' is forced on us by the
fact that handlebodies do not contain annuli which are simultaneously
incompressible and boundary-incompressible.  (Recall that a surface
$F$, properly embedded in a compact manifold $N$, is {\em
  boundary-incompressible} if for every disk $(D, \alpha, \beta)
\subset (N, F, \bdy N)$ such that $\bdy D = \alpha \cup \beta$,
$\alpha \cap \beta$ equals two points, $D \cap F = \alpha$, and $D
\cap \bdy N = \beta$ then there is a disk $(E, \alpha, \beta') \subset
(F, \alpha, \bdy F)$ with $\bdy E = \alpha \cup \beta'$.)
\end{remark}

\begin{remark}
\label{Rem:Mobius}
Suppose that the genus of the handlebody $V$ is greater than one.
Suppose that $A \subset V$ is a properly embedded \Mobius~band.  The
irreducibility of $V$ implies that $\double{A}$ is incompressible.
Also, $\double{A}$ cannot be boundary-parallel as this would imply
that $V$ is a solid torus.  Hence $A$ is essential in $V$.
\end{remark}

If $J \subset V$ is an essential disk, annulus, or \Mobius~band then
$\bdy J \subset \bdy V$ is a collection of essential simple closed
curves.  

\begin{define}
Fix a Heegaard splitting $H \subset M$.  Suppose $J$ and $K$ are each
an essential disk, annulus, or \Mobius~band in the handlebodies $V$ and
$W$, respectively.  If $\bdy J \subset H$ and $\bdy K \subset H$ have
a component in common then $H$ admits a {\em joined $J/K$ pair}.  If
$J$ has a boundary component which is disjoint from some boundary
component of $K$ then $H$ admits a {\em disjoint $J/K$ pair}.
\end{define}

With this terminology a splitting $H$ is {\em reducible} if $H$ admits
a joined disk/disk pair.  If not reducible, then a splitting is {\em
irreducible}.  A formulation of Haken's Lemma~\cite{Haken68} is then:

\begin{lemma}
\label{Lem:Haken}
If $M$ is reducible then every Heegaard splitting $H \subset M$ is
reducible.
\end{lemma}

As a generalization of reducibility Casson and
Gordon~\cite{CassonGordon87} give the following fruitful notion:

\begin{define}
A splitting $H \subset M$ is {\em weakly reducible} if $H$ admits a
disjoint disk/disk pair.  If $H$ is not weakly reducible then it is
{\em strongly irreducible}.
\end{define}

If $H$ admits a joined disk/annulus pair then, after
boundary-compressing the annulus, $H$ admits a disjoint disk/disk
pair.  Thus $H$ is weakly reducible.  If $H$ admits a joined
disk/\Mobius~pair then, as $H$ has genus two or higher, $M$ has
$\RR\PP^3$ as a connect summand and $H$ is reducible.  A somewhat
deeper observation is the ``no nesting'' lemma, due to Scharlemann
(see Lemma~2.2 of~\cite{Scharlemann98}):

\begin{lemma}
\label{Lem:NoNesting}
Suppose that $H \subset M$ is a strongly irreducible splitting and $D
\subset M$ is an embedded disk, with $\interior(D)$ transverse to $H$,
and with $\bdy D \subset H$.  Then $\bdy D$ bounds a disk embedded in
$V$ or $W$.
\end{lemma}

Thompson~\cite{Thompson99} generalizes weak reducibility as follows:

\begin{define}
A Heegaard splitting $H \subset M$ has the {\em disjoint curve
property} if there is an essential simple closed curve $\gamma \subset
H$ and a pair of essential disks $D \subset V$, $E \subset W$ where
$$D \cap \gamma = \gamma \cap E = \emptyset.$$
A splitting is {\em full} if $H$ does not admit such a triple.
\end{define}

\begin{remark}
This notation is prompted by the fact that if $H$ is full then
whenever $D \subset V$ and $E \subset W$ are essential disks the
curves $\bdy D$ and $\bdy E$ {\em fill} the splitting surface.  That
is, $H \setminus (\bdy D \cup \bdy E)$ is a collection of disks.
\end{remark}

\begin{remark}
\label{Rem:Pairs}
If $H$ admits a disjoint disk/annulus pair, a joined annulus/ann\-ulus
pair, or a joined annulus/\Mobius~band pair then $H$ has the disjoint
curve property (DCP).  This is shown by boundary-compressing the
annuli (or the double of the \Mobius~band) to obtain the required
disks.  The required curve $\gamma$ is the appropriate boundary
component of the essential annulus.  (Joined \Mobius~pairs and disjoint
disk/\Mobius~pairs do not arise in our proof.)
\end{remark}

Recall that if $M$ is closed, orientable, irreducible, and contains an
embedded $\pi_1$--injective torus then $M$ is {\em toroidal}.  Work of
Kobayashi~\cite{Kobayashi88a} shows that if $H$ is a strongly
irreducible splitting of a toroidal manifold then $H$ admits a joined
annulus/annulus pair.  Equipped with this Thompson~\cite{Thompson99}
obtains:

\begin{lemma} 
\label{Lem:Thompson}
If $M$ is toroidal then every Heegaard splitting $H \subset M$ has the
disjoint curve property.
\end{lemma}

\begin{figure}[ht!]
\psfrag{R}{R}
\psfrag{W}{W}
\psfrag{D}{D}
\psfrag{I}{I}
\psfrag{S}{S}
\psfrag{F}{F}
$$\begin{array}{c}
\epsfig{file=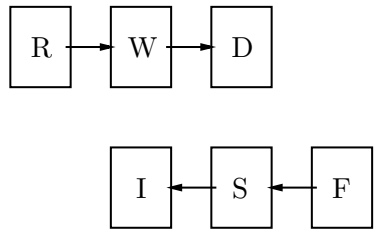, height = 3.0 cm}
\end{array}$$
\caption{Reducibility implies weakly reducibility which in turn
  implies the DCP.  Also fullness implies strong irreducibility which
  in turn implies irreducibility.  It is an amusing exercise to find
  a Heegaard splitting for every column.}
\label{Fig:Implications}
\end{figure}

See Figure~\ref{Fig:Implications} for a diagram of the relations
between the various properties of Heegaard splittings.

We end this section with an observation of Hempel's~\cite{Hempel01}
which relies heavily on the classification of Heegaard splittings of
Seifert fibred spaces~\cite{MoriahSchultens98}:

\begin{lemma}
\label{Lem:Hempel}
If $M$ is a Seifert fibred space then every Heegaard splitting $H
\subset M$ has the disjoint curve property.
\end{lemma}

\section{Normal surface theory}
\label{Sec:Normal}

This section presents the necessary tools from normal surface theory.
For a more complete treatment consult~\cite{JacoOertel84},
\cite{HassEtAl99}, or~\cite{JacoTollefson95}.  See~\cite{Rubinstein97}
for an introduction to almost normal surfaces.

Fix a triangulation, $T_N$, of $N$, a compact orientable
three-manifold.  Denote the $i$--skeleton of $T_N$ by $T_N^i$.  Suppose
$F \subset N$ is an properly embedded surface.  

\begin{define}
The {\em weight} of $F$, $w(F)$, is the number of intersections
between $F$ and the one-skeleton of $T_N$.  The {\em boundary-weight},
$w(\bdy F)$, is the number of intersections between $\bdy F$ and the
one-skeleton of $\bdy N$.
\end{define}

The surface $F$ is {\em normal} if $F$ is transverse to the skeleta of
$T_N$ and intersects each tetrahedron in a collection of {\em normal
  disks}. See Figure~\ref{Fig:NormalDisks} for pictures of the two
kinds of normal disks; the {\em normal triangle} and the {\em normal
  quadrilateral}.  The value of normal surfaces is amply demonstrated
by the following ``normalization'' lemma, essentially due to
Haken~\cite{Haken68}:

\begin{figure}[ht!]
$$\begin{array}{cc}
\epsfig{file=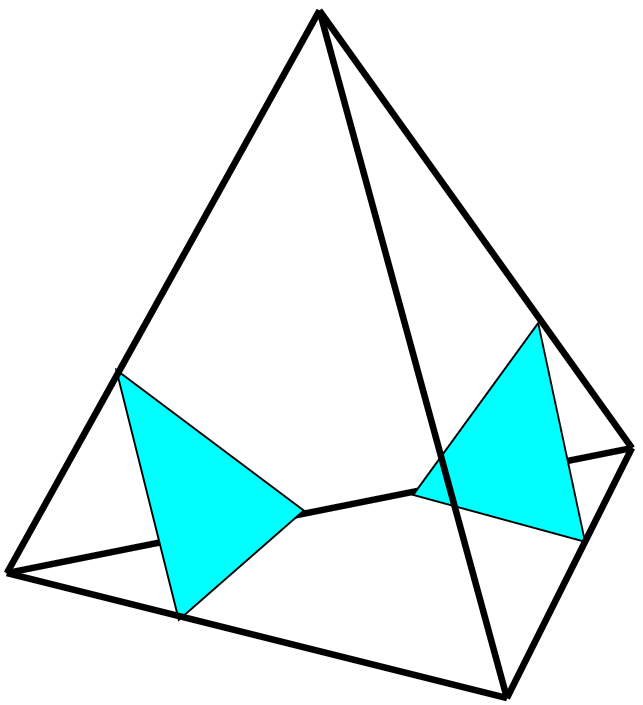, height=3.5cm} &
\epsfig{file=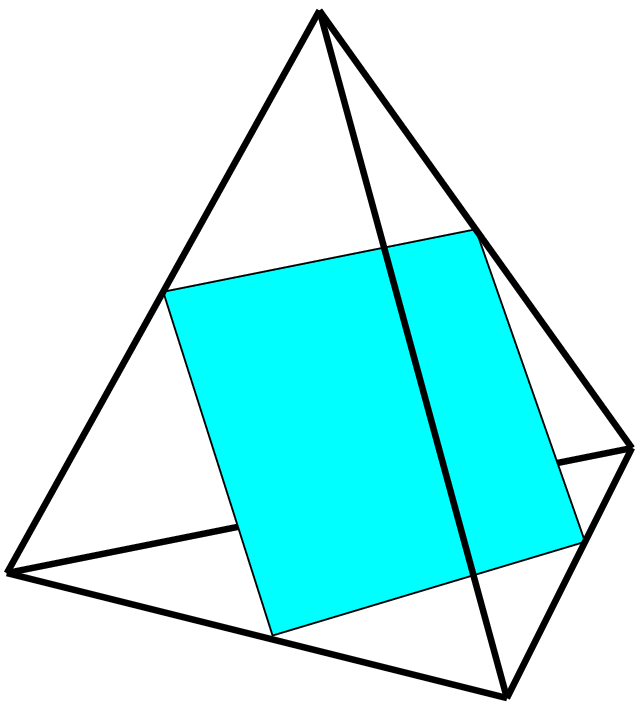, height=3.5cm}
\end{array}$$
\caption{Normal disks}
\label{Fig:NormalDisks}
\end{figure}

\begin{lemma}
\label{Lem:NormalizeIncomp}
Suppose that $(N,T_N)$ is compact, orientable, irreducible, and
boundary-irreducible. If $F \subset N$ is properly embedded,
incompressible, and boundary-incompressible then $F$ is properly
isotopic to a normal surface $F'$ with $w(F') \leq w(F)$ and with
$w(\bdy F') \leq w(\bdy F)$.  Furthermore, if $F \cap T^2_N$ contains
an arc with both endpoints in a single edge of some face, then also
have $w(F') < w(F)$.
\end{lemma}

\begin{figure}[ht!]\small
\psfrag{F}{$F$}
$$\begin{array}{c}
\epsfig{file=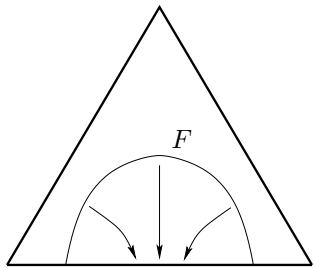, height = 3.5 cm}
\end{array}$$
\caption{A weight-reducing isotopy}
\label{Fig:BentArcWithIsotopy}
\end{figure}

(Recall that a three-manifold $N$ is {\em boundary-irreducible} if
$\bdy N$ is incompressible in $N$.)  To illustrate the final sentence
of the lemma see Figure~\ref{Fig:BentArcWithIsotopy}.  Note that the
``normalization'' lemma holds even for one-sided $F$.

A {\em normal isotopy} of $N$ fixes every simplex of $T_N$ setwise.
Two normal surfaces are {\em equivalent} if there is a normal isotopy
taking one to the other.  Two normal surfaces, $J, K \subset N$, are
{\em compatible} if in each tetrahedron of $T_N^3$ the surfaces $J$
and $K$ have the same types of quad (or one or both have no quads).
In this situation form the {\em Haken sum} $F = J + K$ as follows:

Normally isotope $J$ to make $J$ transverse to $K$, to make $\Gamma =
J \cap K$ transverse to the skeleta of $T_N$, and to minimize the
number of components of $\Gamma$.  The curves $\gamma \subset \Gamma$
are the {\em exchange curves}.  For every such curve $\gamma$ let
$R(\gamma)$ be the closure (taken in $N$) of $\neigh_N(\gamma)$, an
open regular neighborhood of $\gamma$.  Since $N$ is assumed
orientable $R(\gamma)$ is a solid torus or a one-handle, with $\gamma$
as its core curve.

\begin{figure}[ht!]\small
\psfrag{Ar}{$A_r$}
\psfrag{Ai}{$A_i$}
$$\begin{array}{c}
\epsfig{file=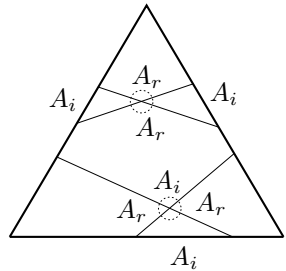, height=3.5cm}
\end{array}$$
\caption{Regular and irregular annuli}
\label{Fig:SmallNeighborhood}
\end{figure}

\begin{figure}[ht!]\small
\psfrag{F}{$J$}
\psfrag{G}{$K$}
\psfrag{Arg}{$A_r(\gamma)$}
\psfrag{gamma}{$\gamma$}
\psfrag{regularexchange}{Cut-and-paste}
\psfrag{exchangeband}{$C(\gamma)$}
$$\begin{array}{c}
\epsfig{file=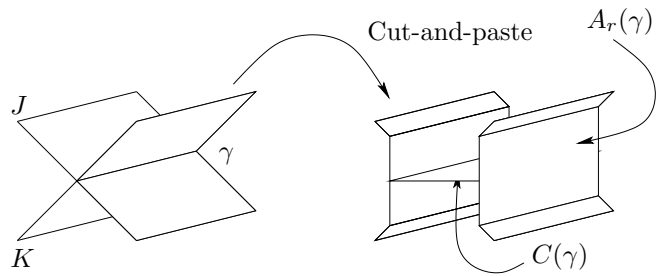, height=3.5cm}
\end{array}$$
\caption{The exchange band for $\gamma$}
\label{Fig:ExchangeBand}
\end{figure}

Now, $\frontier_N(R(\gamma)) \setminus (J \cup K)$, the frontier of
$R(\gamma)$ in $N$ minus the two surfaces, is a union of annuli or
rectangles because $N$ is orientable. Taking closures, divide these
into two sets, the {\em regular bands} $A_r(\gamma)$ and the {\em
  irregular bands} $A_i(\gamma)$, as indicated by
Figure~\ref{Fig:SmallNeighborhood}.  (The fact that $J$ and $K$ are
compatible proves that each band receives the same label from every
face of $T^2_N$.)  Finally, as in Figure~\ref{Fig:ExchangeBand}, form the
Haken sum:
$$F = \left( (J \cup K) \setminus \bigcup_\Gamma \{\neigh_N(\gamma)\} \right)
             \cup \bigcup_\Gamma \{A_r(\gamma)\}.$$

\begin{define}
Each connected component of $(J \cup K) \setminus \bigcup_\Gamma
\{\neigh_N(\gamma)\}$ is a {\em patch} of the sum $J + K$ while each
regular band $A_r(\gamma)$ is a {\em seam}. 
\end{define}

Note that $F$ is again a normal surface which, up to normal isotopy,
does not depend on the choices permitted by the above construction.
Also, the weight, boundary-weight, and Euler characteristic of $F$ are
additive with respect to Haken sum.

\begin{remark}
\label{Rem:Seams}
If $A$ is a seam of $F = J + K$ then $A$ is not a \Mobius~band.
Equivalently, if $A \cap \bdy N = \emptyset$ and $\alpha$ is the core
curve of $A$ then $\alpha$ preserves orientation in $F$.  This is
trivial when, as with us, $N$ is orientable.  (However it also holds
for $N$ non-orientable.)
\end{remark}

A normal surface is {\em fundamental} if it admits no such
decomposition.  A crucial finiteness result due to Haken (and recorded
in~\cite{JacoTollefson95}) is:

\begin{lemma}
\label{Lem:FiniteFundamental}
For a fixed compact, triangulated three-manifold $(N^3, T_N)$ there are
only finitely many fundamental normal surfaces, up to normal isotopy.
\end{lemma}

A bit of linear algebra and a few estimates give a more quantitative
result (see Hass \etal~\cite{HassEtAl99}):

\begin{lemma}
\label{Lem:WeightOfFundamental}
If $F \subset (N, T_N)$ is fundamental then $w(F) \leq \exp(14|T_N^3|)$.
Also, $\euler(F) \geq -\exp(14|T_N^3|)$. 
\end{lemma}

Here $\exp(x) = 2^x$ is the exponential function and $|T_N^3|$ is the
number of tetrahedra in the triangulation $T_N$.  Abusing notation
$|\cdot|$, below, also denotes the number of components of a
topological space.

Each ``cut-and-paste'' operation involved in the Haken sum $F = J + K$
may be recorded by an embedded {\em exchange band} $C(\gamma)$. The
band $C(\gamma)$ is defined as follows: Choose a homeomorphism between
$R(\gamma)$, the neighborhood of the intersection curve $\gamma$, and
a square-bundle over $\gamma$.  This square-bundle has in each fibre
(\ie~each square) $J \cup K$ as the two diagonals, $A_r(\gamma)$ as
one set of opposing sides, and $A_i(\gamma)$ as the other set of
sides.  The exchange band $C(\gamma)$ is the surface which intersects
each square in the line segment connecting midpoints of the
$A_r(\gamma)$ sides.  See Figure~\ref{Fig:ExchangeBand}.

Each $C(\gamma)$ is an annulus (\Mobius~band) when $\gamma$ is closed
and orientation-preserving (reversing) in $J$ and $K$.  The band
$C(\gamma)$ is a rectangle when $\gamma$ is not closed.  Every seam
$A_r(\gamma)$ is a closed regular neighborhood (in $F$) of a boundary
component of some $C(\gamma)$.

A sum $F = J + K$ is {\em reduced} if the normal surface $F$ cannot be
realized as a sum $J' + K'$ where $J'$ and $K'$ are again normal,
properly isotopic to $J$ and $K$ respectively, with $|J' \cap K'| < |J
\cap K|$. Note that the isotopy between $J$ and $J'$ ($K$ and $K'$)
need not be normal.  The ``no disk patches'' lemma is a version of a
key technical result proved by Jaco and Oertel (Lemma~2.1
of~\cite{JacoOertel84} ):

\begin{lemma}
\label{Lem:NoDiskPatches}
Suppose that $(N, T_N)$ is a triangulated, compact, orientable,
irreducible, boundary-irreducible three-manifold.  Suppose that $F
\subset (N, T_N)$ is a properly embedded, incompressible,
boundary-incompressible normal surface which minimizes the
lexicographic complexity $(w(\bdy F), w(F))$ in its isotopy class.
Suppose $F = J + K$ is a reduced sum.  Then no patch of $J + K$ is a
disk.
\end{lemma}


Notice that $F$ need not be two-sided.  Also, this paper only uses the
lemma in the situation $\bdy K = \emptyset$.

\begin{remark}
\label{Rem:ExchangeOnTunnel}
Also contained in \cite{JacoOertel84} is the following observation.
Suppose $F = J + K$ is given, as above.  Then no exchange annulus $A$
of $J + K$ will be parallel into $F$.  For suppose that $(A, \bdy A)
\subset (N, F)$ was parallel through a solid torus, in the complement
of $F$, to an annulus $(C, \bdy C) \subset (F, \bdy A)$.  Then let
$F''$ be a slight push-off of $(F \setminus C) \cup A$ (see
Figure~\ref{Fig:ExchangeOnTunnel}) and note that either $w(F'') <
w(F)$ or the final line of Lemma~\ref{Lem:NormalizeIncomp} applies to
$F''$.  In either case a contradiction is obtained.
\end{remark}

\begin{figure}[ht!]
$$\begin{array}{cc}
\epsfig{file=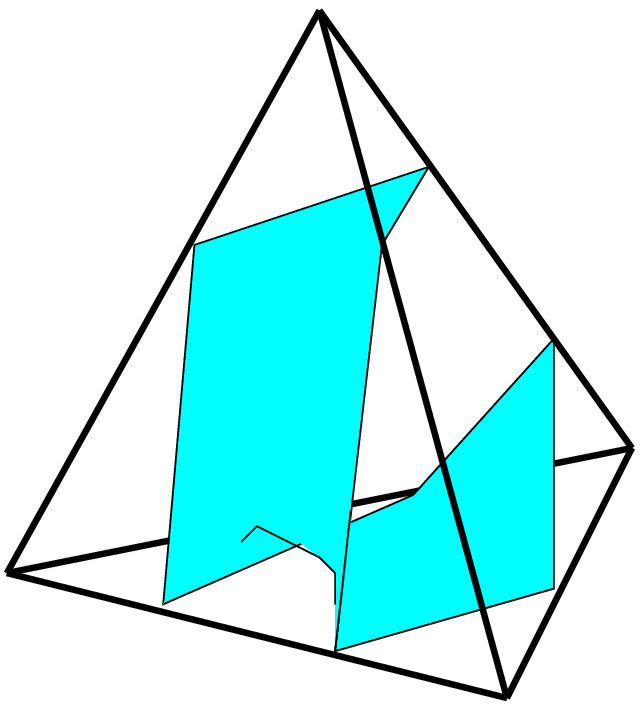, height=3.5cm} &
\epsfig{file=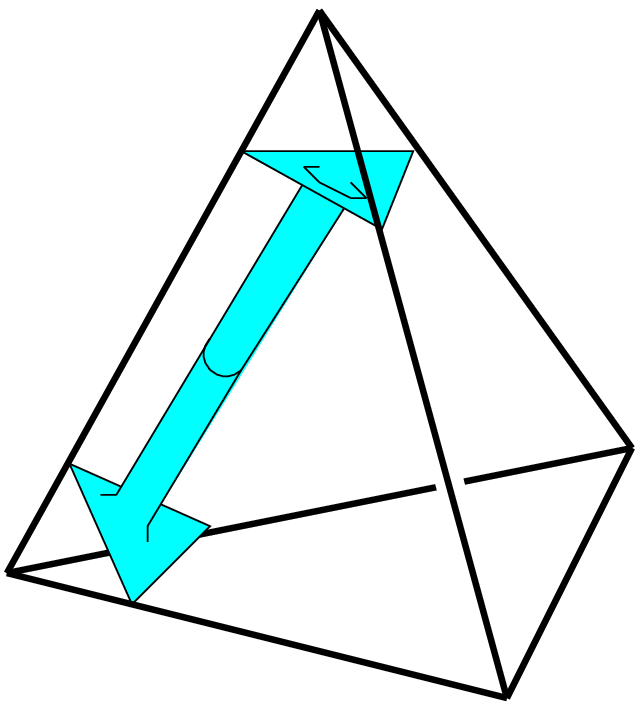, height=3.5cm}
\end{array}$$
\caption{Almost normal pieces}
\label{Fig:AlmostNormalPieces}
\end{figure}

A surface $H$, embedded in a closed manifold $(M, T_M)$, is {\em almost
normal} if $H$ is transverse to the skeleta of $T_M$ and intersects
each tetrahedron but one in a collection of normal disks. In the
remaining tetrahedron $H$ yields a collection of normal disks and at
most one {\em almost normal piece}. Two of the five kinds of almost
normal pieces are shown in Figure~\ref{Fig:AlmostNormalPieces}.  The first
is the {\em almost normal octagon} while the second is one of the
{\em almost normal annuli}.  The other possibilities are an almost
normal annulus between triangles of the same type, between a pair of
quads, or between a triangle and a quad.

The following theorem (see the papers of
Rubinstein~\cite{Rubinstein97} and Stocking~\cite{Stocking00}) is
another important normalization result:

\begin{theorem}
\label{Thm:NormalizeStrIrred}
If $H \subset (M,T_M)$ is a strongly irreducible Heegaard splitting
then $H$ is isotopic to an almost normal surface.
\end{theorem}

\section{Blocked and shrunken submanifolds}
\label{Sec:Blocks}

This section introduces blocks and deals with submanifolds (of a
triangulated three-manifold) which are naturally contained in unions
of blocks.  

Suppose that $\tau$ is a regular Euclidean tetrahedron. Let $S \subset
\tau$ be a disjoint collection of normal disks together with at most
one almost normal piece. ($S = \emptyset$ is allowed.)

\begin{define}
A {\em block}, $B$, is the closure of a connected component of $\tau
\setminus S$.  Suppose $B$ is adjacent to exactly two normal disks,
$D$ and $E$.  Suppose also that $D$ is normally isotopic in $\tau$
to $E$. Then $B$ is a {\em product block}. All other blocks are called
{\em core blocks}.
\end{define}

\begin{figure}[ht!]\small
$$\begin{array}{cc}
\epsfig{file=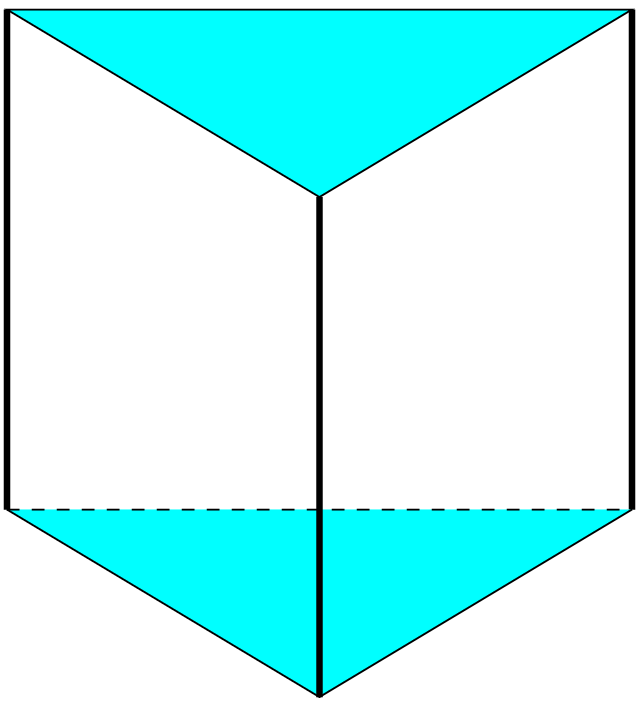, height=3.3cm} &
\epsfig{file=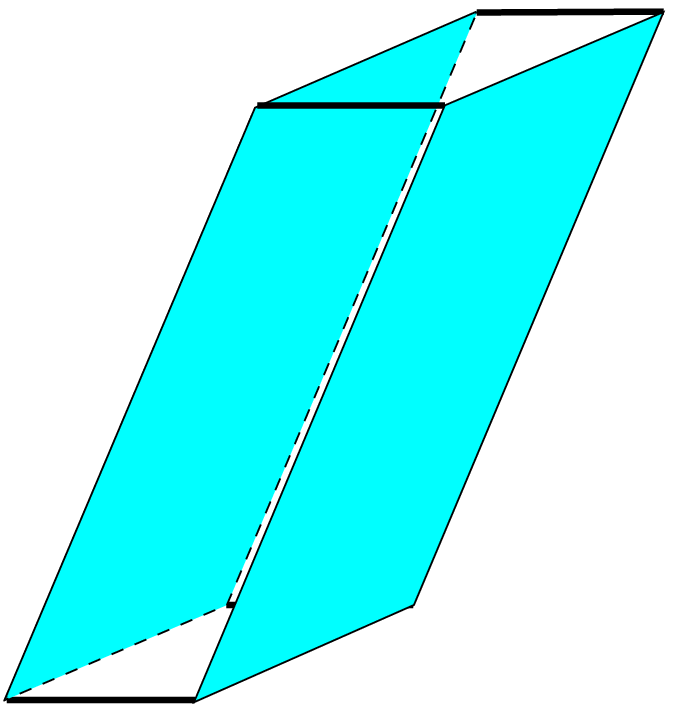, height=3.3cm} \\
\end{array}$$
\begin{center}
The product blocks
\end{center}
$$\begin{array}{cc}
\epsfig{file=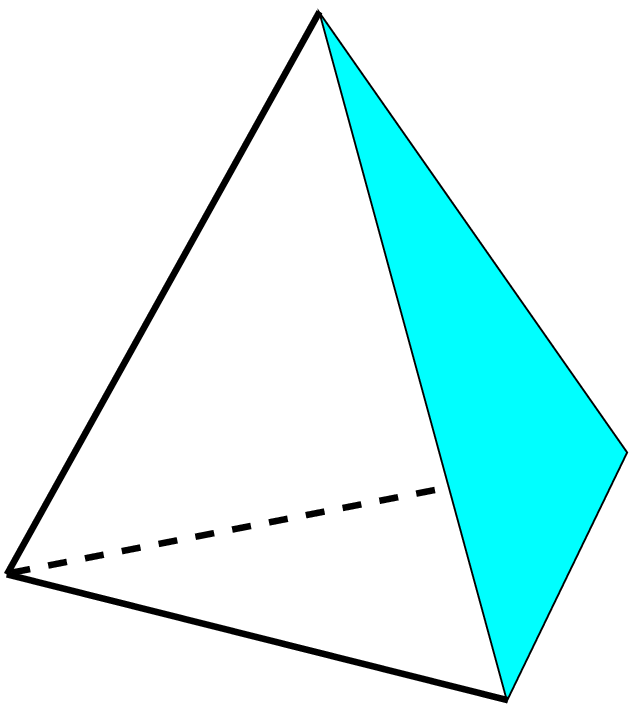, height=3.3cm} &
\epsfig{file=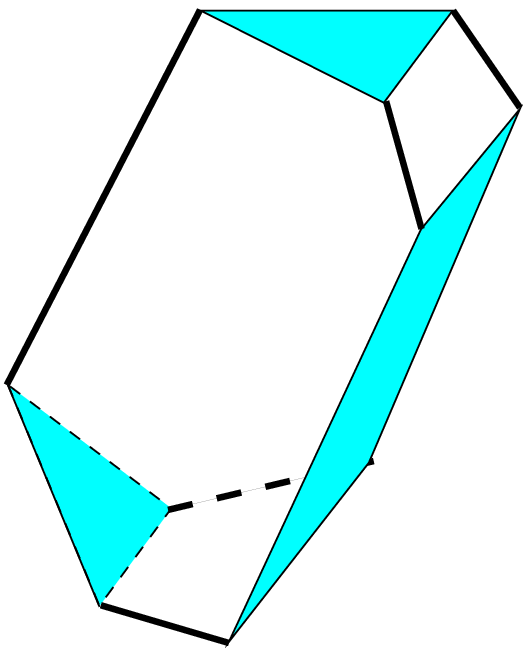, height=3.3cm} \\
\epsfig{file=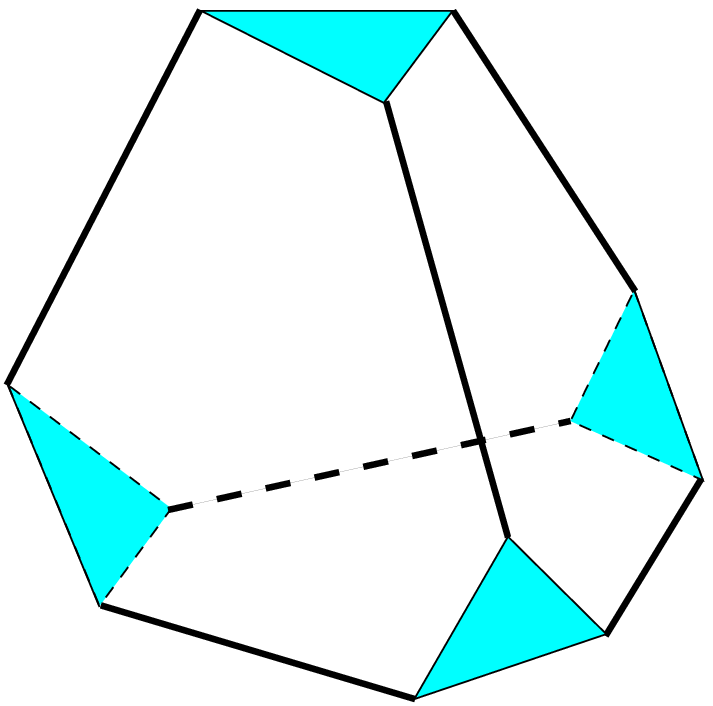, height=3.3cm} &
\epsfig{file=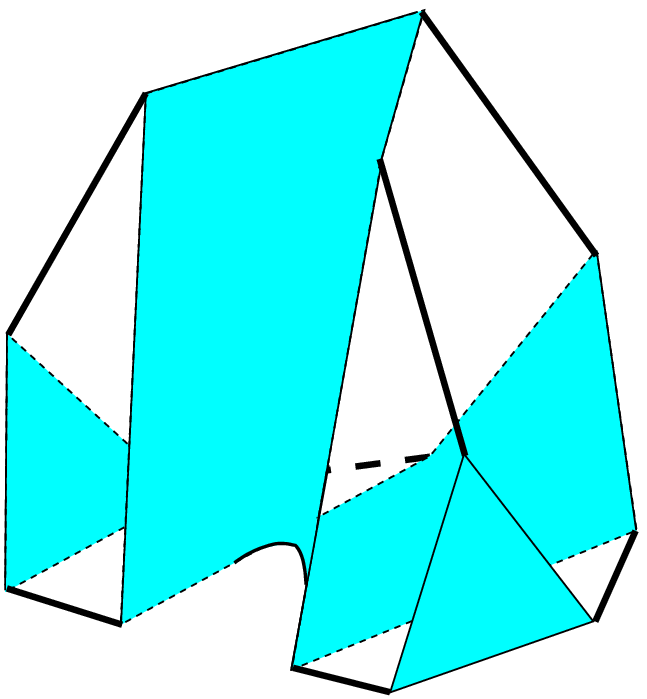, height=3.3cm} \\
\epsfig{file=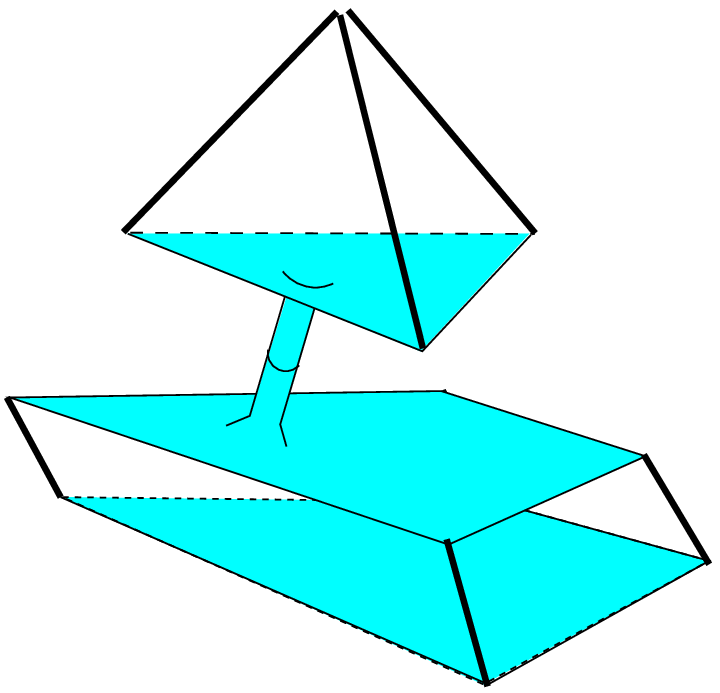, height=3.3cm} &
\epsfig{file=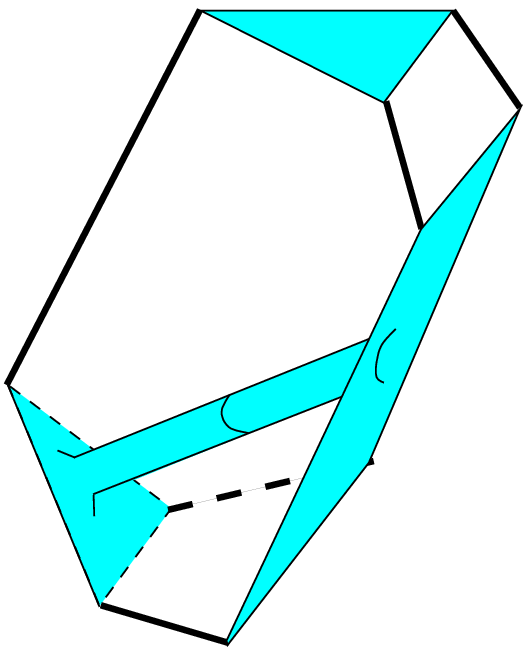, height=3.3cm} \\
\end{array}$$
\begin{center}
A few core blocks
\end{center}
\caption{Some of the possible blocks (not all to scale)}
\label{Fig:Blocks}
\end{figure}

See Figure~\ref{Fig:Blocks} for pictures of the two product blocks and
six of the many possible core blocks. Note that a few of the faces are
shaded --- these are the faces which lie in the surface $S$.  Such
faces form the {\em horizontal boundary} of the block while the faces
$B \cap \tau^2$, in the two-skeleton, are the {\em vertical boundary}
of $B$.  Define a {\em vertical rectangle} of a block $B$ to be any
face of the vertical boundary which is disjoint from the zero-skeleton
of $T_M$ and which has four sides.  For example, all vertical faces of
a product block are rectangles.

Let $(M,T_M)$ be a closed, orientable, triangulated three-manifold.

\begin{define}
A three-dimensional submanifold $X \subset (M,T_M)$ is {\em blocked} if
$\bdy X$ is an almost normal surface and $X$ is a union of
blocks.
\end{define}

As a specific example, if $H \subset (M,T_M)$ is connected, separating,
normal surface then the closure of a component of $M \setminus H$ is a
blocked submanifold.

\begin{remark}
\label{Rem:BoundedNumberOfCoreBlocks}
Note that a blocked submanifold contains at most $6 t + 1$ core
blocks, where $t = |T_M^3|$ is the number of tetrahedra in $T_M$.  To see
this, note that a tetrahedron contains at most five parallel families
of normal disks. Thus any tetrahedron not meeting an almost normal piece
contributes at most six core blocks. The almost normal piece (if it
exists) gives the final core block.
\end{remark}


Suppose that $X \subset (M,T_M)$ is blocked.  Suppose $Y \subset X$ is
a subset of $X$.  Define the {\em horizontal boundary} of $Y$ to be
$\bdy_h Y = Y \cap \bdy X$.  Also define the {\em vertical
  boundary} of $Y$ to be $\bdy_v Y = \fr_X(Y)$, the set of points in
$Y$ which do not have an open neighborhood (relative to $X$) in $Y$.

Recall that $\neigh_X(\cdot)$ denotes a open regular neighborhood
taken inside of $X$.  

\begin{define}
Let $X \subset (M,T_M)$ be blocked. A three-dimensional submanifold $Y
\subset X$ is {\em shrunken} if there is a union of blocks $\blocky{Y}
\subset X$, with
$$Y = \blocky{Y} \setminus \neigh_X(\bdy_v \blocky{Y}).$$
\end{define}

Thus the shrunken $Y$ is obtained by removing a small neighborhood of
$\bdy_v \blocky{Y}$ from $\blocky{Y}$.  (That is, by {\em shrinking}
$\blocky{Y}$.) Note that $Y$ uniquely determines $\blocky{Y}$ as well
as the reverse.

\begin{define}
If $Y$ is any submanifold of a blocked $X \subset (M, T_M)$ then
define the {\em complexity} of $Y$:
$$c(Y) = \euler(Y) + |\bdy_v Y|.$$
\end{define}

This is an exact parallel of the {\em complexity} of a compact surface
$G$, $c(G) = \euler(G) + |\bdy G|$.  If $Y$ is an $I$--bundle
then $c(Y)$ agrees with the complexity of the base surface. 

If $Y$ is connected and $c(Y) \leq 0$ then $Y$ is {\em
large}. Again, this parallels the notation for surfaces: A connected
surface $G$ is large if $G$ has nonpositive complexity.
Equivalently, $G$ is large if $G$ contains an orientation-preserving,
nonseparating, simple closed curve.

We end this section with a discussion of how to obtain new shrunken
submanifolds from old. Let $Y$ be a shrunken submanifold of the
blocked $X$. Let $B \subset X$ be a block which non-trivially
intersects the interior of $Y$. Then $Y \setminus \neigh_X(B)$ is
again a shrunken submanifold inside of $X$.  The following lemma
records how the complexity of $Y \setminus \neigh_X(B)$ differs from
that of $Y$:

\begin{lemma}
\label{Lem:IBundleComplexityChange}
Given a blocked $X \subset (M,T_M)$, a shrunken submanifold $Y \subset
X$, and a block $B$ of $X$ which non-trivially intersects $Y$ then:
$$c(Y \setminus \neigh_X(B)) \leq c(Y) + \blockconst.$$
\end{lemma}


\begin{proof}
First note that the block $B$ contains at most four vertices of
$T_M^0$.  Also, $B \cap (T_M^1 \setminus T_M^0)$ is the union of at
most $12$ connected arcs.  This is because $B$ intersects any single
edge of $T_M^2$ at most twice.  Finally, $\bdy_v B$ contains at most
$12$ faces.  To see this last, note that the dual graph to $B \cap
T_M^2$ has as at most $12$ edges, dual to the components of $B \cap
(T_M^1 \setminus T_M^0)$, and this graph has no vertices of valency
zero or one.

To obtain $Y \setminus \neigh_X(B)$ remove neighborhoods of
certain skeleta of $B$ from $Y$.  To be precise, for every vertex in
$B \cap T_M^0 \cap Y$ remove a small open neighborhood of the vertex
from $Y$.  Let $Y'$ be the compact submanifold obtained.  Each vertex
removed increases the Euler characteristic by one and contributes one
additional vertical boundary component.  Thus,
$$c(Y') \leq c(Y) + 8.$$
Remove from $Y'$ a regular neighborhood of $(B \cap Y') \cap
T_M^1$, giving the submanifold $Y''$.  Note that $c(Y'') \leq c(Y')$ as
any increase in $|\bdy_v|$ is exactly compensated by a decrease in
Euler characteristic.

Finally remove a regular neighborhood of $(B \cap Y'') \cap T_M^2$ to
obtain $Y'''$.  There are at most $12$ faces to remove; each increases
Euler characteristic and $|\bdy_v|$ by at most one.  Thus,
$$c(Y''') \leq c(Y'') + 24.$$
Deleting the interior of $B$ yields $Y \setminus \neigh_X(B)$ and does
not increase the complexity.  To see this note that the Euler
characteristic remains unchanged when $B$ is a solid torus and
actually decreases when $B$ is a ball.  Also, removing the interior of
$B$ discards at least one component of the vertical boundary.  Thus
$c(Y \setminus \neigh_X(B)) \leq c(Y) + \blockconst$, as desired.
\end{proof}

\section{The main proposition}
\label{Sec:MainProposition}

Here the main technical proposition is given.  As the proof of the
proposition is quite long it is broken apart among the sections
which follow.  

\begin{proposition}
\label{Prop:AlmostNormalPlusHighGenusImpliesDCP}
Suppose $(M,T_M)$ is a closed, orientable, triangulated three-manifold.
Suppose $H \subset M$ is an almost normal Heegaard splitting with
genus $g(H) > \dcpconst$.  Then $H$ has the disjoint curve
property.
\end{proposition}


Here $t = |T_M^3|$ is the number of tetrahedra in $T_M$ and, again,
$\exp(x) = 2^x$ is the exponential function.

\begin{remark}
The estimate on genus is extremely wasteful.  However, our methods
cannot be used to give upper bounds better than $\exp(c\cdot t^2)$
where $c$ is a constant.  See also the discussion in
Section~\ref{Sec:Conjecture}.
\end{remark}

We now outline the proof of
Proposition~\ref{Prop:AlmostNormalPlusHighGenusImpliesDCP}, which operates
via contradiction.

First assume that an almost normal Heegaard splitting $H$ is given
which both has large enough genus and is full.  Both hypotheses are
used throughout the proof.  Cut $M$ along $H$ to obtain a pair of
blocked handlebodies, $V$ and $W$.  Section~\ref{Sec:Shelling} and
Section~\ref{Sec:CappingOff} produce a shrunken $I$--bundle $Y \subset
V$ with $\bdy_h Y \subset H$.  The $I$--bundle $Y$ has $\euler(Y) \leq
-\exp(2^{15} t^2)$ while still having a small number of fairly short
vertical boundary annuli.  Furthermore, it is shown that the vertical
and horizontal boundaries of $Y$ are incompressible inside of $V$.

\begin{remark}
\label{Rem:YUntwisted}
Note that, when $Y$ is not a twisted $I$--bundle, the proof becomes
quite straight-forward.  It is easy to find an essential annulus in
$W$ with one boundary component in $\bdy_h Y$.  Simply take this curve
and ``push it through'' $Y$ to obtain the desired joined
annulus/annulus pair.
\end{remark}

Next, in Section~\ref{Sec:FAndN}, a manifold $N \homeo M \setminus
\neigh_M(\bdy_v Y)$ is introduced and shown to be irreducible and
boundary-irreducible.  Also, a section $F$ of the $I$--bundle $Y$ is
chosen and shown to be incompressible and boundary-incompressible
inside of $N$.

The manifold $N$ is triangulated in Section~\ref{Sec:NormalizingF}.  The
derived triangulation $T_N$ has at most $2^{9} t^2$ tetrahedra where
$t = |T_M^3|$.  This bound follows from the small size of $\bdy_v Y$.
The surface $F$ is then normalized with respect to $T_N$ so that
$w(\bdy F) \leq 2^{6} t^2$.  This, together with the fact that
$\euler(F) = \euler(Y)$, allows a Haken decomposition of $F$ of a
pleasant form.  Note that $F$ cannot be a fundamental normal surface
in $N$, as the Euler characteristic of $F$ is too negative.

Finally, in Section~\ref{Sec:LastAnnulus}, the exchange bands of the
decomposition of $F$ are examined.  It is one of these which gives
rise to an joined annulus/annulus (or \Mobius/annulus) pair for the
splitting surface $H$.  Boundary-compressing each of these essential
annuli (or the double of the \Mobius~band) shows that $H$ had the
disjoint curve property, a final contradiction. This will complete the
proof of Proposition~\ref{Prop:AlmostNormalPlusHighGenusImpliesDCP}.

\section{Shelling high genus splittings}
\label{Sec:Shelling}

This section deals with one of the main objects of interest: almost
normal Heegaard splittings.  A pair of blocked handlebodies and a few
shrunken submanifolds will be derived from such a splitting. These
submanifolds will be large $I$--bundles and their vertical annuli
provide our first supply of essential annuli.  The constructions and
notations here introduced will be used heavily throughout the proof of
Proposition~\ref{Prop:AlmostNormalPlusHighGenusImpliesDCP}.

Suppose that $(M,T_M)$ is a closed, orientable, triangulated
three-manifold and $H$ is an almost normal splitting with genus
greater than $\dcpconst$ as given by
Proposition~\ref{Prop:AlmostNormalPlusHighGenusImpliesDCP}. Recall
that $t = |T^3_M|$.  Suppose also that $H$ has the smallest weight
among all almost normal surfaces isotopic to $H$.  Finally suppose,
for an eventual contradiction, that $H$ is a full splitting.

First analyze the basic structure of $M$.  From Haken's Lemma
(Lemma~\ref{Lem:Haken} above) and the fact that full splittings are
irreducible it follows that:

\begin{claim}
\label{Clm:MIrreducible}
The manifold $M$ is irreducible.
\end{claim}

From Waldhausen's Theorem~\cite{Waldhausen68} further deduce:

\begin{claim}
\label{Clm:MNotS3}
The manifold $M$ is not homeomorphic to the three-sphere.
\end{claim}

Recall that $H$ divides $M$ into two blocked handlebodies, $V$ and
$W$.  Let $\blocky{V_1}$ be the union of all product blocks in $V$.
Let $\blocky{W_2}$ be the union of all product blocks $B \subset W$
such that $\bdy_h B \subset \bdy_h \blocky{V_1}$.  Let $V_1$ and $W_2$
be the corresponding shrunken submanifolds inside of $V$ and $W$.  As
$M$ is orientable, $V_1$ and $W_2$ are orientable submanifolds.  Also,
by construction, $\bdy_h W_2 \subset \bdy_h V_1$.

Let $Z_2 \subset W_2$ be a connected component which has smaller (more
negative) complexity than any other component of $W_2$.  Amongst
the one or two components of $V_1$ meeting $\bdy_h Z_2$, let $Y_1$ be
one with smallest complexity.

We investigate the properties of these $I$--bundles, in several
steps.

\begin{claim}
The $I$--bundle $W_2$ has $c(W_2) \leq 1 - g(H) + 20 \cdot \blockconst
t$.  Also, $|\bdy_v W_2| < 104 t$.  
\end{claim}

\begin{proof}
The $I$--bundle $W_2$ is obtained from $W$ by removing all core blocks
from $W$ and then removing all product blocks of $W$ which are
adjacent to core blocks of $V$, across a horizontal face.  There are
at most $6t + 1$ of the former by
Remark~\ref{Rem:BoundedNumberOfCoreBlocks}.  It remains to estimate the
number of the latter.  If the almost normal piece sits inside the
tetrahedron $\tau$ then $\tau$ contains at most six families of
parallel product blocks.  If not, then $\tau$ contains at most five
such families.  Each family of parallel product blocks contains at
most two blocks adjacent to a core block.  Thus, the $I$--bundle $W_2$
is obtained by removing at most another $10 t + 2$ product blocks from
$W$.  Applying Lemma~\ref{Lem:IBundleComplexityChange} the complexity of
$W_2$ is bounded by
$$c(W_2) \leq 1 - g(H) + 20 \cdot \blockconst t.$$
To obtain the second inequality of the claim: note that each core
block has at most $8$ vertical rectangles: those four-sided vertical
faces which do not meet the zero-skeleton of $T_M$.  Product blocks
have either $3$ or $4$ vertical rectangles.  Each contributes at most
one to the ``length'' of $\bdy_v W_2$, so it certainly follows that 
$|\bdy_v W_2| \leq 8(6t + 1) + 4(10t + 2) \leq 104 t$.
\end{proof}


Recall that $Z_2$ was chosen to be a component of $W_2$ with smallest
(most negative) complexity.

\begin{remark}
By the above claim the $I$--bundle $W_2$ contains at most
$104 t$ components.  Thus, as complexity is additive over
disjoint union, the component $Z_2$ has complexity 
$$c(Z_2) \leq \frac{c(W_2)}{104 t} 
         \leq \frac{1 - g(H)}{104 t} + 8.$$
\end{remark}

Recall that $Y_1$ is a component of $V_1$ such that $\bdy_h Y_1$
contains a component of $\bdy_h Z_2$. 

\begin{remark}
The complexity of $Y_1$ is bounded above by:
$$c(Y_1) \leq \frac{1}{2} c(Z_2) 
         \leq \frac{1 - g(H)}{208 t} + 4.$$ 
This follows from Lemma~\ref{Lem:SurfaceComplexityChange}, proved in the
appendix, applied to the relevant components of $\bdy_h Y_1$ and
$\bdy_h Z_2$.  Conclude that $c(Y_1) \leq -\exp(2^{15} t^2)$. 
\end{remark}

The $I$--bundles $Y_1$ and $Z_2$ both contain vertical nonseparating
annuli, say $A \subset Y_1$ and $B \subset Z_2$.  Applying
Lemma~\ref{Lem:NonSepVerticalAnnuliAreEssential} from the appendix find:

\begin{claim}
\label{Clm:AnnuliForY1}
There are essential annuli $A \subset V$, $B \subset W$ such that
$\bdy A$ and one component of $\bdy B$ are embedded in $\bdy_h Y_1$.
\end{claim}

From this immediately deduce a non-essential but technically useful
corollary.

\begin{claim}
\label{Clm:NoAlmostNormalAnnulus}
The almost normal splitting $H$ does not contain an almost normal
annulus.
\end{claim}


\begin{proof}
Recall that $H$ was assumed, at the beginning of this section, to be
the least weight surface among all almost normal surfaces isotopic to
$H$.  Suppose that the splitting $H$ contains an almost normal
annulus.  Then there is a disk $D \subset M$ such that $\bdy D$ is a
core curve for the almost normal annulus, $D$ meets $H$ only along
$\bdy D$, and the disk $D$ does not meet the two-skeleton $T_M^2$.

There are two possibilities: either $\bdy D$ is inessential in $H$ or
it is essential.  If inessential then surger $H$ along $D$ to obtain
$H'$. The surface $H'$ is the union of a normal two-sphere $S$ and a
normal surface $H''$.  As $M$ is irreducible
(Claim~\ref{Clm:MIrreducible}) the two-sphere $S$ bounds a three-ball $U
\subset M$.  If the splitting $H$ were contained in $U$ then $M \homeo
S^3$.  This is impossible by Claim~\ref{Clm:MNotS3}.  Deduce that
there is an isotopy taking $H''$ to $H$ which is supported on a
regular neighborhood of $U$.  But this is another impossibility as
$H''$ is normal with $w(H'') < w(H)$.

It follows that $\bdy D$ is essential in $H$.  Note that $\bdy D
\cap \bdy_h V_1 = \emptyset$.  By Claim~\ref{Clm:AnnuliForY1}, $D$ forms
one-half of a disjoint disk/annulus pair.  Boundary-compressing the
given annulus demonstrates that $H$ has the DCP, another
contradiction.
\end{proof}


This proof exemplifies a theme of this paper --- always seek to
compare a troublesome disk to a vertical annulus provided by one of
the $I$--bundles.  As $H$ is full it cannot admit a disjoint
disk/annulus pair.  Conclude that the vexatious disk does not exist.

\section{Capping off to obtain $Y$}
\label{Sec:CappingOff}

This section discusses the annuli forming the vertical boundary of
$Y_1$.  Some of these bound two-handles in $V$.  These two-handles
will be added to $Y_1$ to form a bigger $I$--bundle, $Y$.

For every vertical annulus $C_i \subset \bdy_v Y_1$, with both
components of $\bdy C$ inessential in $H$, let $D_i$ and $D'_i$ be the
disks in $H$ which the components of $\bdy C_i$ bound. The two-sphere
$C_i \cup D_i \cup D'_i$ bounds a ball, $U_i \subset V$, as
handlebodies are irreducible.  Note that $U_i \cap Y_1 = C_i$.  (If
not $Y_1 \subset U_i$ and $\bdy_h Y_1$ would be planar, contradicting
the fact that $Y_1$ is large.)

Give each ball $U_i$ the structure of an $I$--bundle over the disk
$\DD^2$. Arrange matters so that the $I$--bundle structure on $C_i$
coming from $Y_1$ and $U_i$ agree. Let $Y = Y_1 \cup \bigcup_i
U_i$. Again, $Y$ is a shrunken submanifold of $V$, $Y$ is a large
$I$--bundle, and $\bdy_h Y_1 \subset \bdy_h Y$.  We are
now equipped to make several observations about the structure of $Y$.

First of all, $c(Y) = c(Y_1) \leq -\exp(2^{15} t^2)$.  Also, $\bdy_v
Y$ is non-empty.  (Because surface groups are not free!)

\begin{remark}
\label{Rem:AnnuliForY}
Note that $\bdy_h Y_1 \subset \bdy_h Y$.  Thus
Claim~\ref{Clm:AnnuliForY1} applies to $Y$ just as it applies to
$Y_1$.  That is, there are essential annuli $A \subset V$, $B \subset
W$ such that $\bdy A$ and one component of $\bdy B$ are embedded in
$\bdy_h Y$. 
\end{remark}

\begin{define}
The {\em corners} of the $I$--bundle $Y$ are the components of $\bdy
(\bdy_v Y)$.  
\end{define}

\begin{claim}
\label{Clm:CornersOfY}
No corner of $Y$ bounds an embedded disk in $M$ which is transverse to
$H$.  (So, in particular, the corners of $Y$ are essential in $H$.)
\end{claim}

\begin{proof}
Fix attention on $C$, a component of $\bdy_v Y$.  Let $\bdy C =
\gamma \cup \gamma'$.  Suppose, for a contradiction, that $\gamma$
bounds a disk in $M$ which is transverse to $H$.  By taking the union
of this disk with $C$ and performing a small isotopy $\gamma'$ also
bounds a disk in $M$, transverse to $H$.

By the construction of $Y$ one component of $\bdy C$ is essential in
$H$.  Apply ``no nesting'' (Lemma~\ref{Lem:NoNesting}) to this
essential boundary component.  It follows that $Y$ has a corner which
is essential in $H$ and bounds a disk in $V$ or $W$.  By
Remark~\ref{Rem:AnnuliForY} there is a disjoint disk/annulus pair and so
$H$ has the DCP.  This is a contradiction.
\end{proof}


It immediately follows that:

\begin{claim}
\label{Clm:VerticalBdyIncompressibleInV}
The components of $\bdy_v Y$ are incompressible in $V$.
\end{claim}

Also:

\begin{claim}
\label{Clm:HorizontalBdyIncompressibleInV}
The components of $\bdy_h Y$ are incompressible in $V$.
\end{claim}

\begin{proof}
Let $D \subset V$ be a compressing disk for $\bdy_h Y$ which
minimizes $|D \cap \bdy_v Y|$. Then $D$ does not meet $\bdy_v Y$ by
Claim~\ref{Clm:VerticalBdyIncompressibleInV}. Thus $D$ lies inside of
$Y$. But this is impossible, as the horizontal boundary of an
$I$--bundle is always incompressible inside the $I$--bundle.
\end{proof}

\section{Focusing on the middle surface}
\label{Sec:FAndN}

This section discusses the geometric structure of a section $F \subset
Y$.  This surface is properly embedded in the submanifold $N =
\closure{M \setminus R}$, as discussed below.

\begin{figure}[ht!]\small
\psfrag{Y}{$Y$}
\psfrag{V-YunionR}{$V \setminus (Y \cup R)$}
$$\begin{array}{c}
\epsfig{file=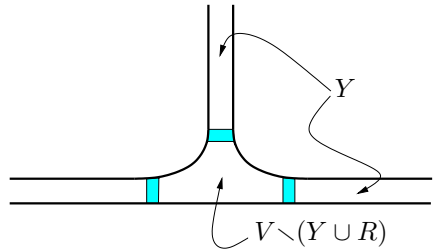, height = 3.3 cm}
\end{array}$$
\caption{A cross-section of $V$.  The union of tori, $R$, is shaded.}
\label{Fig:YInV}
\end{figure}

Let $Y \subset V$ be the $I$--bundle produced above, in
Section~\ref{Sec:CappingOff}.  Let $R$ be the following union of closed
solid tori:
$$R = \closure{\neigh_{\closure{V \setminus Y}}(\bdy_v Y)}.$$ 
Then $R \cap Y = \bdy_v Y$.  That is, if $C$ is a component of $\bdy_v
Y$ and $R'$ is the component of $R$ meeting $C$ then $R'$ is a collar
for $C$ inside of $\closure{V \setminus Y}$.  See Figure~\ref{Fig:YInV}.

\begin{figure}[ht!]\small
\psfrag{R}{$R$}
\psfrag{BdyR}{$\bdy R$}
\psfrag{F}{$F$}
\psfrag{BdyY}{$\bdy_h Y$}
$$\begin{array}{c}
\epsfig{file=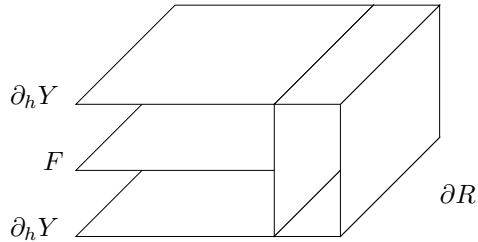, height = 3.2 cm}
\end{array}$$
\caption{Side-view of $F \subset Y$}
\label{Fig:SliceOfYAndR}
\end{figure}

Let $F \subset Y$ be a properly embedded surface, transverse to the
fibres of $Y$, meeting each fibre exactly once, and with $F \cap
\bdy_h Y = \emptyset$.  See Figure~\ref{Fig:SliceOfYAndR}.  Let $\pi\from
Y \to F$ be the natural bundle map: $\pi^{-1}(x)$ is the $I$--fibre of
$Y$ meeting the point $x \in F$.  The surface $F$ is a section of the
$I$--bundle $Y$.  

It follows that the complexity of $F$ equals that of $Y$.  The Euler
characteristic of $F$ is less than $c(F) = c(Y) \leq -\exp(2^{15}
t^2)$ because $\bdy F$ is nonempty.  Note that $F$ is properly
embedded in the three-manifold $N = \closure{M \setminus R}$.  When
discussing $Y$ inside of $N$ we will regard it as a closed regular
neighborhood of $F$.  The rest of this section studies the pair $F
\subset N$.

\begin{claim}
\label{Clm:FIncompressible}
The surface $F$ is incompressible in $N$.
\end{claim}

Before proving this, observe that a stronger statement would be that
$\bdy_h Y$ is incompressible in $W$. However, this does not seem to
necessarily hold when $Y$ is a twisted $I$--bundle.

\begin{proof}[Proof of Claim~\ref{Clm:FIncompressible}]
Suppose that $F$ is compressible, via a disk $D \subset N$.  Let
$\alpha = \bdy D$ and set $A = \pi^{-1}(\alpha)$.  Note that $A$ is an
annulus, not a \Mobius~band.  We may assume, after a small isotopy of
$D$ $\rel~\alpha$, that $D \cap Y = \neigh_D(\alpha) \subset A$.

Let $D' = D \setminus \interior(Y)$.  Set $\alpha' = \bdy D' \subset
\bdy_h Y$.  The curve $\alpha'$ is essential in $H$.  (If $\alpha'$
bounds a disk in $\bdy_h Y$ then $\alpha$ is null-homotopic in $F$, a
contradiction.  If $\alpha'$ bounds a disk in $H$ meeting $\bdy
(\bdy_h Y)$ then have a contradiction to Claim~\ref{Clm:CornersOfY}.)

Lemma~\ref{Lem:NoNesting} replaces $D'$ by a disk $D''$ with $\bdy D''
= \bdy D' = \alpha'$ and with the interior of $D''$ disjoint from the
splitting $H$.  Note that $D''$ cannot lie inside of $V$ by
Claim~\ref{Clm:HorizontalBdyIncompressibleInV}.  It follows that $D''$
lies in $W$.

Now choose $\beta \subset F \setminus \alpha$ a simple closed curve
which is nonseparating and orientation-preserving in $F$.  (Such a
curve exists as $c(F)$ is very negative.)  By
Lemma~\ref{Lem:NonSepVerticalAnnuliAreEssential} the annulus $B =
\pi^{-1}(\beta)$ is essential in $V$.  Thus $D''/B$ forms a disjoint
disk/annulus pair and $H$ has the DCP.  This is a contradiction.
\end{proof}


\begin{claim}
\label{Clm:FBdyIncompressible}
The surface $F$ is boundary-incompressible in $N$.
\end{claim}


\begin{proof}
Suppose that $F$ is boundary-compressible, via a disk $D \subset N$.
Choose $D$ to minimize the quantity $|\interior(D) \cap H|$.  Let
$\alpha = (\bdy D) \cap F$.  We may assume, after a small proper
isotopy of $D$ $\rel~\alpha$ which leaves $|\interior(D) \cap H|$
unchanged, that $D \cap Y = \neigh_D(\alpha) \subset
\pi^{-1}(\alpha)$.  Let $D' = D \setminus \interior(Y)$.  Set $\alpha'
= D' \cap \bdy_h Y$.

Note that $D' \cap H$ contains no simple closed curves embedded in
$\interior(D')$.  For suppose to the contrary that $\delta$ is such an
innermost such curve and $\delta$ bounds a subdisk $E \subset D'$.
Now $\delta$ must be essential in $H$.  (If not $\delta$ bounds a disk
disk $E' \subset H$.  As $\bdy_h Y$ is nonplanar $E' \cap \bdy_h Y =
\emptyset$.  So use an innermost disk of $E'$ to surger $D$, reducing
$|\interior(D) \cap H|$.  This is a contradiction.)  Thus $E$ is an
essential disk in $V$ or $W$ which is disjoint from $\bdy_h Y$.  By
Remark~\ref{Rem:AnnuliForY} the splitting $H$ admits a disjoint
disk/annulus pair, implying that $H$ has the DCP, a contradiction.
So, as desired, $D' \cap H$ contains no simple closed curves embedded
in $\interior(D')$.

Consider the components of $D' \setminus H$.  The closure of
one of these contains $\alpha'$.  Denote the closure of that component
by $D''$.  See Figure~\ref{Fig:DAndD''}.

\begin{figure}[ht!]\small
\psfrag{D''}{$D''$}
\psfrag{a}{$\alpha$}
\psfrag{a'}{$\alpha'$}
$$\begin{array}{c}
\epsfig{file=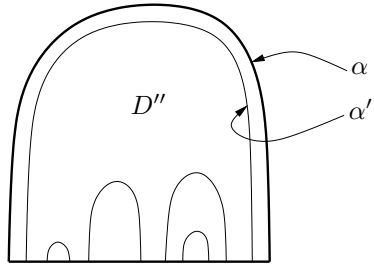, height = 3.5 cm}
\end{array}$$
\caption{The disk $D$, showing $D \cap H$ and $D''$}
\label{Fig:DAndD''}
\end{figure}

Note that $D''$ is a disk embedded in $W$.  Let $\beta \subset F$ be a
non-separating, orientation-preserving simple closed curve which is
disjoint from $\alpha$.  (Such a curve exists as $c(F)$ is very
negative.)  Set $B = \pi^{-1}(\beta)$; the annulus which is the union
of all fibres of $Y$ meeting $\beta$.  By
Lemma~\ref{Lem:NonSepVerticalAnnuliAreEssential} the annulus $B$ is
essential in $V$.  

If $D''$ were essential then $D''$ and $B$ would form a disjoint
disk/annulus pair, a contradiction.  Thus $\bdy D''$ bounds a disk $E
\subset H$ with $\bdy E \cap Y = \alpha'$.  As $\bdy(\bdy_v Y) \cap E$
contains no simple closed curves (Claim~\ref{Clm:CornersOfY}), it
follows that $\alpha'$ is boundary-parallel in $\bdy_h Y$.  So
$\alpha$ is homotopic (rel boundary) into $\bdy F$.  This implies that
the original disk $D$ was not a boundary-compression after all.
\end{proof}


Before examining the three-manifold $N$ further we will need the
following:

\begin{claim}
\label{Clm:YNotInBall}
The $I$--bundle $Y \subset M$ is not contained inside of a three-ball
in $M$.
\end{claim}

\begin{proof}
Suppose that $Y$ is contained inside of a three-ball.  Pick a
three-ball $U \subset M$ such that $U$ contains $Y$ and minimizes $|S
\cap H|$, where $S = \bdy U$.  If $|S \cap H| = 0$ then $H$ is
contained inside of a ball in $M$.  Since handlebodies are irreducible
this implies that $M$ is the three-sphere, contradicting
Claim~\ref{Clm:MNotS3}.

Thus $S \cap H$ is nonempty.  Let $\delta \subset (S \cap H)$ be
innermost in $S$.  Let $D \subset S$ be an innermost disk which
$\delta$ bounds.  If $\delta$ is essential in $H$ then, by
Remark~\ref{Rem:AnnuliForY}, $D$ is half of a disjoint disk/annulus
pair and $H$ has the DCP, a contradiction.

Thus $\delta$ bounds a disk in $H$.  Let $E \subset H$ be an innermost
disk of $S \cap H$.  Set $\epsilon = \bdy E$.  The disk $E$ is
disjoint from $Y$, as $\bdy_h Y$ is nonplanar.  Let $E'$ and $E''$ be
slight push-offs of $E$.  Perform disk surgery on $S$ along $E$ to
obtain $S' \disjoint S'' = (S \setminus \neigh_S(\epsilon)) \cup E'
\cup E''$, a disjoint union of two-spheres.  Note that both $|S' \cap
H|$ and $|S'' \cap H|$ are less than $|S \cap H|$.  Also, as $M$ is
irreducible, both $S$ and $S'$ bound three-balls.

If $E \subset U$ then both $S'$ and $S''$ bound balls inside of
$U$ and one of these contains $Y$. This however contradicts the
minimality of $|S \cap H|$.  Suppose instead that $E \cap U = \epsilon
= \bdy E$.  If one of $S'$ or $S''$ bounds a ball containing $U$
then that one of $S'$ or $S''$ bounds a ball containing $Y$, a
contradiction to minimality.  It follows that both bound balls
disjoint from $U$.  Thus $M$ is homeomorphic to the three-sphere.
This last contradicts Claim~\ref{Clm:MNotS3}.
\end{proof}


From this, deduce:

\begin{claim}
\label{Clm:NIrreducible}
The three-manifold $N$ is irreducible.
\end{claim}

\begin{proof}
Suppose not.  Let $S \subset N$ be a reducing sphere for $N$ which
minimizes the quantity $|S \cap F|$.  

Suppose that $|S \cap F| = 0$.  After a small isotopy of $Y$ we may
assume that $S \cap Y = \emptyset$, because $Y$ is a regular
neighborhood of $F$ in $N$.  Now, $S$ bounds a three-ball $U \subset
M$, as $M$ is irreducible (Claim~\ref{Clm:MIrreducible}).  By
Claim~\ref{Clm:YNotInBall} the $I$--bundle $Y$ does not meet $U$.  As
$R$ is a small collar of $\bdy_v Y$ deduce that $U$ is disjoint from
the union of solid tori $R$.  Thus $U$ embeds in $N = \closure{M
  \setminus R}$ and $S$ was not a reducing sphere, a contradiction.

Suppose instead that $|S \cap F| > 0$.  Let $\delta \subset (S \cap
F)$ be innermost in $S$, bounding the disk $D \subset S$. As $F$ is
incompressible (Claim~\ref{Clm:FIncompressible}) the curve $\delta$
also bounds a disk in $F$.  Let $E \subset F$ be an innermost subdisk
with boundary $\bdy E = \epsilon \subset S \cap F$. As in the proof of
Claim~\ref{Clm:YNotInBall} above, surger $S$ along the disk $E$ to
obtain two-spheres $S'$ and $S''$.

If neither of these is a reducing sphere then a standard argument
shows that $S$ was also not a reducing sphere, a contradiction.  But
$S'$, say, being a reducing sphere contradicts the minimality
assumption.  
\end{proof}

Finally, deduce that $N$ is boundary-irreducible, or what amounts
to the same thing:

\begin{claim}
\label{Clm:BdyNIsIncompressible}
The boundary of $N$ is incompressible inside of $N$.
\end{claim}

\begin{proof}
Suppose that $\bdy N$ compresses.  Let $D$ be a compressing disk for
$\bdy N$ which minimizes $|D \cap F|$.  Note that $D \cap F$ contains
no simple closed curves.  (If such exists then by the
incompressibility of $F$ there is a disk $E \subset F$ with $E \cap D
= \bdy E$.  Surger $D$ along $E$ to reduce $|D \cap F|$.)

Suppose first that $|D \cap F| = 0$.  Properly isotope $D$ in $N$ to a
disk $D'$ so that $\bdy D'$ lies inside some component of $R \cap H$.
That is, $\bdy D'$ lies in $H$ and is parallel to a component of $\bdy
(\bdy_h Y)$.  This contradicts Claim~\ref{Clm:CornersOfY}.

Suppose that $|D \cap F| > 0$.  Let $D'' \subset D$ be an outermost
bigon with respect to the arcs $D \cap F$.  As $F$ is
boundary-incompressible the arc $D'' \cap F$ cuts $F$ into two pieces,
one of which is a bigon, $E$.  An outermost bigon $E' \subset E$ gives
a $\bdy$--surgery for $D$.  Performing this surgery will produce at
least one compressing disk for $\bdy N$ which meets $F$ in fewer arcs,
a contradiction.
\end{proof}

\section{Normalizing $F$}
\label{Sec:NormalizingF}

This section finds a triangulation for $N$ and shows that $F$ may be
normalized inside of $N$.  This normal surface will then decompose as
a well-behaved Haken sum.  As $F$ is not necessarily two-sided, it is
not immediately obvious how to replace the technicalities of
triangulations and normal surface theory with a neater proof using
branched surfaces.  However, I would be happy to be corrected on this
point.


\begin{claim}
\label{Clm:TriangulationOfN}
The manifold $N$ may be triangulated with fewer than $2^{9} t^2$
tetrahedra so that the surface $F$ has boundary-weight less than
$2^{6} t^2$.
\end{claim}


Here $t = |T_M^3|$ is the number of tetrahedra in the given
triangulation of $M$.  The proof operates by constructing a
sufficiently simple triangulation, $T_N$, of $N$.  This construction
is complicated while also being somewhat uninteresting.  I would
advise the reader to take the above claim on faith in his or her first
reading of this section and instead skip directly to
Section~\ref{Sec:MakeFNormal}.

\subsection{Constructing the Graph $\Theta_F$}  
Let $\blocky{Y}$ be the union of blocks meeting $Y$.  So $Y \subset
\blocky{Y}$.  As usual, define $\bdy_h \blocky{Y} = \blocky{Y} \cap
H$.  The subset $\blocky{Y}$ is not a submanifold of $M$ but it does
have the structure of an $I$--bundle.  So let $\blocky{F}$ be a section
of $\blocky{Y}$ which contains $F$ and which is disjoint from $\bdy_h
\blocky{Y}$.

\begin{figure}[ht!]
\psfrag{}{}
\psfrag{}{}
$$\begin{array}{c}
\epsfig{file=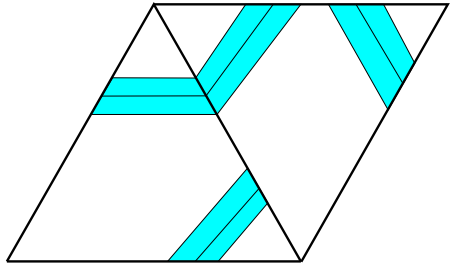, height = 3.5 cm}
\end{array}$$
\caption{A bit of $\Theta_F$ with the nearby $\Theta_Y$ shaded}
\label{Fig:GraphAndNeighborhood}
\end{figure}

Define $\Theta_Y = \bdy_v \blocky{Y} =
\closure{\frontier_M(\blocky{Y}) \setminus \bdy_h \blocky{Y}}$.  This
is a regular neighborhood, taken in $T_M^2$, of the graph $\Theta_F =
\frontier_M(\blocky{F})$ which is also contained in $T_M^2$.  (See
Figure~\ref{Fig:GraphAndNeighborhood}.)

Before continuing on we will tidy up the graph $\Theta_F$ a bit.
Recall that every tetrahedron of $T_M$ is modelled on the regular
Euclidean tetrahedron of side-length one.  Apply a normal isotopy of
$M$ (\ie~an isotopy fixing every simplex of $T_M$ setwise) to make
every normal disk of $\blocky{F}$ meeting $\Theta_F$ {\em straight}.
That is, after the isotopy, every normal triangle of $\blocky{F}$
which meets $\Theta_F$ is a geodesic triangle in its containing
tetrahedron while every quad is either a geodesic quadrilateral or a
union of two geodesic triangles.  We also require that every edge of
$\Theta_F$ lies outside of a $1/3$--neighborhood of the vertices of
$T_M^0$.

Our next goal will be to bound the complexity of the graph $\Theta_F$.

\subsection{Bounding the Complexity of $\Theta_F$}  
\begin{proofclaim}
$\Theta_F$ has at most $4 t$ edges and at most $4 t$ vertices.
\end{proofclaim}

\begin{proof}
Every edge of $\Theta_F$ lies inside a vertical rectangle of
$\Theta_Y$ and each of these is simultaneously a vertical rectangle of
some product block of $\blocky{Y}$ and of some core block of $V$.
Thus to bound the edges of $\Theta_F$ it suffices to bound the number
of vertical rectangles appearing in the boundaries of core blocks.

Recall that the splitting $H$ does not contain an almost normal
annulus, by Claim~\ref{Clm:NoAlmostNormalAnnulus}.  Thus each core
block in $V$, disjoint from the almost normal octagon, meets
$\blocky{Y}$ along at most two vertical rectangular faces.  The core
block which contains the almost normal octagon, if it exists, meets
$\blocky{Y}$ in at most four rectangular faces.  Also, at most $2 t$
of the core blocks in $V$ have a vertical rectangle as a face.  (See
Figure~\ref{Fig:Blocks}.)

Now, the tetrahedron containing the octagon contains only one core
block with vertical rectangles.  Deduce that $\Theta_Y = \bdy_v
\blocky{Y}$ is a union of at most $4 t$ rectangles lying in the
two-skeleton.  It follows that the number of edges of $\Theta_F$ is at
most $4t$.  We now must bound the number of vertices in $\Theta_F$.

Consider a vertex $v$ of $\Theta_F^0$.  As $\Theta_F$ is a union of
edges, the vertex $v$ has valence one or higher.  Suppose $v$ has
valence one exactly.  Now, $v$ sits on some edge, $e$, of the
triangulation $T_M^1$ and $e$ intersects exactly one vertical
rectangle, $E$, of $\Theta_Y$.  Thus in every tetrahedron containing
$e$ there is a product block of $\blocky{Y}$ in that tetrahedron
containing $v$.  It follows that $E$ intersects the interior of $Y$, a
contradiction.  It follows that $v$ must have valence two or higher.
Thus $\Theta_F$ has at most $4t$ vertices.
\end{proof}

\subsection{Triangulating $\bdy R$}  
We will now triangulate $\bdy R = \bdy N$ as well as a ``crust'',
$\Delta_F$, which connects $\Theta_F$ to a collection of longitudes
for $R$.  This will permit an efficient triangulation of $N$.

\begin{figure}[ht!]\small
\psfrag{l}{$\lambda_i$}
\psfrag{BdyR}{$\bdy R$}
\psfrag{R}{$R$}
\psfrag{D_Y}{$\Delta_F$}
$$\begin{array}{c}
\epsfig{file=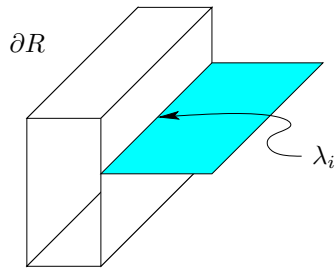, height = 3.5 cm}
\end{array}$$
\caption{A cross-sectional view of $\bdy R$ ($\Delta_F$ shaded)}
\label{Fig:NearDelta}
\end{figure}

Let $\Delta_F = \closure{\blocky{F} \setminus (R \cup F)}$.  Then the
frontier of $\Delta_F$ is the union of $\Theta_F$ and a collection of
longitudes, $\{ \lambda_i \}$, for the tori $R_i \subset R$.  See
Figure~\ref{Fig:NearDelta}.

Choose $\epsilon > 0$ very small so that $\neigh_N(\Theta_F, 5
\epsilon)$, the $5 \epsilon$ neighborhood of $\Theta_F$, is a regular
neighborhood, and so that every edge of $\Theta_F$ has length at least
$100\epsilon$.  Perform another normal isotopy of $M$, fixing
$\blocky{F}$ setwise, so that the longitudes $\lambda_i = \Delta_F
\cap \bdy R_i$ lie on the boundary of the handlebody
$\neigh_N(\Theta_F, \epsilon)$.

Partition each longitude $\lambda_i$ into two sets of arcs as follows.
First for every open tetrahedron $\tau_j \in T_M^3$ there is a
collection of arcs $\{ \alpha_i^{j,k} \} = \lambda_i \cap
\neigh(\tau^1_j, 2\epsilon)$.  Here $\tau^1_j$ is the one-skeleton of
the tetrahedron $\tau_j$ while $j$ runs from $1$ to $|T_M^3|$.  The
arcs $\{ \beta_i^{j,k} \}$ are the components of $\closure{\lambda_i
\setminus (\cup_{j,k} \alpha_i^{j,k}})$.  Count the number of
$\beta$'s and $\alpha$'s:

\begin{proofclaim}
The number of $\beta$'s is at most $4 t$ while the number of
$\alpha$'s is at most $24 t^2$.
\end{proofclaim}

\begin{proof}
For every edge $e \in \Theta_F$ there is exactly one of the
$\beta_i^{j,k}$'s lying in the closed $\epsilon$--neighborhood of
$e$.  This bounds the number of $\beta$'s.

Fix $v \in \Theta_F^0$.  Suppose that $v \in e$ where now $e$ is an
edge of $T_M^1$, the one-skeleton of $M$.  Note that $e$ appears as an
edge of triangles in $T_M^2$ in at most $6t$ ways (this is a very
crude bound).  So there are at most $6t - 1$ product blocks in
$\blocky{Y}$ containing the vertex $v$.  It follows that there are at
most $6t - 1$ of the $\alpha_i^{j,k}$'s within the
$2\epsilon$--neighborhood of $v$.  It follows that there are at most
$4t \cdot 6t$ of the $\alpha$'s in total.
\end{proof}

Perform a final normal isotopy fixing $H \cup T_M^2$ pointwise, fixing
the points $\bdy \alpha_i^{j,k}$ for all $i, j, k$, and sending the
arcs $\{ \alpha_i^{j,k} \} \cup \{ \beta_i^{j,k} \}$ to Euclidean line
segments.

\begin{figure}[ht!]\small
\psfrag{BdyR}{$\bdy R$}
\psfrag{D_Y}{$\Delta_F$}
\psfrag{G_Y}{$\Theta_F$}
\psfrag{Alpha}{$\alpha_i^{j,k}$}
\psfrag{Beta}{$\beta_i^{j,k}$}
$$\begin{array}{c}
\epsfig{file=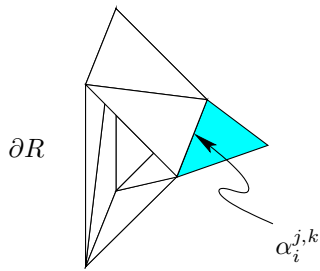, height = 3.5 cm}
\end{array}$$
\caption{Triangulation of $\bdy R$ near $\alpha_i^{j,k}$ ($\Delta_F$ shaded)}
\label{Fig:NearAlpha}
\end{figure}

\begin{figure}[ht!]\small
\psfrag{BdyR}{$\bdy R$}
\psfrag{D_Y}{$\Delta_F$}
\psfrag{G_Y}{$\Theta_F$}
\psfrag{Alpha}{$\alpha_i^{j,k}$}
\psfrag{Beta}{$\beta_i^{j,k}$}
$$\begin{array}{c}
\epsfig{file=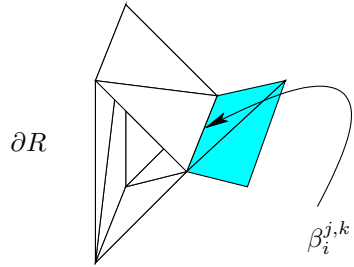, height = 3.5 cm}
\end{array}$$
\caption{Triangulation of $\bdy R$ near $\beta_i^{j,k}$ ($\Delta_F$ shaded)}
\label{Fig:NearBeta}
\end{figure}

Triangulate $\bdy R_i$ and $\Delta_F$ near $\alpha_i^{j,k}$ as shown in
Figure~\ref{Fig:NearAlpha}.  Triangulate $\bdy R_i$ and $\Delta_F$
near $\beta_i^{j,k}$ as shown in Figure~\ref{Fig:NearBeta}.  It follows
that:

\begin{proofclaim}
The triangulation of $\bdy R$ contains at most $6(4 t + 24 t^2)$
triangles.  Also, the triangulation of $\Delta_F$ contains at most $8
t + 24 t^2$ triangles.  Finally, $\bdy F$ meets the one-skeleton of
$\bdy R$ in at most $2(4 t + 24 t^2)$ points. 
\end{proofclaim}

So the boundary-weight of $F$ in $N$ will be at most $2(4 t + 24 t^2)
< 2^{6} t^2$, as desired.

\subsection{Triangulating $N$}  
We will now subdivide the balls $\{ \tau \setminus (R \cup \Delta_F)
\}$ to obtain a triangulation of $N$.

\begin{figure}[ht!]
$$\begin{array}{c}
\epsfig{file=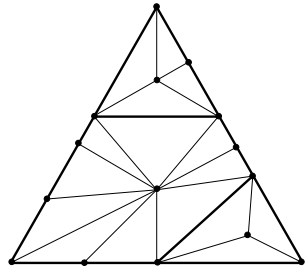, height = 3.5 cm}
\end{array}$$
\caption{The subdivision of a face $\sigma$}
\label{Fig:ConingOverSigmaPrime}
\end{figure}

Every vertex of $T_M^0$ and $\Theta_F^0$ is included in the vertex set
of $T_N$. Fix a face $\sigma \subset T_M^2$.  Let $\delta_1, \delta_2,
\delta_3$ be the edges of $\sigma$.  Let $\sigma'$ be a component of
$\sigma \setminus \Theta_F$.  Then, for each $\sigma'$, add a point of
the interior of $\sigma'$ to the vertex set of $T_N$.  Also, cone from
this point to the boundary of $\sigma'$.  See
Figure~\ref{Fig:ConingOverSigmaPrime}.  Thus $\sigma$ is subdivided
into the following number of triangles:
$$2|\sigma \cap \interior(\Theta_F^1)| + 3 + \sum_{1 = 1}^3 |\delta_i
\cap \Theta_F^0|.$$
For each tetrahedron $\tau \subset T_M^3$ let $\tau'$ be the closure
of $\tau \setminus (R \cup \neigh(\Delta_F))$.  Then $\tau'$ is
homeomorphic to a round ball and the chosen triangulation of $T_M^2
\cup \Delta_F \cup \bdy R$ induces a triangulation on the boundary of
$\tau'$.  To finish, choose a point in the interior of $\tau'$ and
cone from this point to the boundary of $\tau'$.  Doing this for every
tetrahedron $\tau \subset T_M^3$ gives the desired triangulation,
$T_N$.  Note that each triangle of the subdivided $\sigma \subset
T_M^2$ gives us two tetrahedra in $T_N$ as does every triangle in the
$\Delta_F$. However, each triangle of $\bdy R$ gives us only one
tetrahedron in $T_N$.  Adding, find that $|T_N^3|$ is at most
$$
2 \sum_\sigma \left( 2|\sigma \cap \interior(\Theta_F^1)| 
                     + 3 + \sum_{1 = 1}^3 |\delta_i \cap \Theta_F^0| \right) 
+ 2(8 t + 24 t^2) 
+ 6(4 t + 24 t^2).
$$
The above quantity is at most:
\begin{gather*}
16t + 12t + 2(6t \cdot 4t) + 16 t + 48 t^2 + 24 t + 144 t^2\\
= 68 t + 240 t^2 < 2^9 t^2. 
\end{gather*}
Here $t = |T_M^3|$ is the number of tetrahedra in $T_M$.  Deduce that
$N$ admits a triangulation, $T_N$, with $|T_N^3| < 2^{9} t^2$ where
$F$ has boundary-weight $w(\bdy F) < 2^{6} t^2$.  This completes the
proof of Claim~\ref{Clm:TriangulationOfN}.

\subsection{Normalizing the middle surface, $F$}
\label{Sec:MakeFNormal}

\begin{claim}
\label{Clm:MakeFNormal}
There is a proper isotopy $\cf : F \times I \rightarrow N$ such that
$F' = \cf_1(F)$ is normal with respect to $T_N$, $w(F') \leq w(F)$,
$w(\bdy F') \leq w(\bdy F)$.
\end{claim}

All weights here are measured with respect to the triangulation $T_N$.

\begin{proof}
Recall that $F$ is incompressible and boundary-incompressible and that
$N$ is irreducible and boundary-irreducible.  Thus Haken's
normalization procedure (Lemma~\ref{Lem:NormalizeIncomp} in
this paper) gives the desired isotopy.
\end{proof}

Make the further assumption that $F'$, as given by
Claim~\ref{Clm:MakeFNormal}, minimizes the lexicographic complexity
$(w(\bdy F'), w(F'))$.  Abusing notation slightly again use $F$
to denote the normal surface $F'$.  Decompose $F$ as a sum of
fundamental surfaces of the triangulation $T_N$:
$$F = \sum_i m_i G_i + \sum_i n_i F_i + \sum_i p_i T_i.$$
Here the
$G_i$ are fundamental surfaces with nonempty boundary, the $F_i$ are
closed fundamental surfaces with Euler characteristic negative, and
the $T_i$ are also closed, fundamental, but with Euler characteristic
non-negative.  Note that, by Claim~\ref{Clm:TriangulationOfN} and the
lexicographic minimality of $F$, the sum $\sum m_i$ is bounded above
by $2^{6} t^2$. This is because each $G_i$ has non-trivial boundary
and boundary-weight is additive under Haken sum.

Euler characteristic is also additive under Haken sum.  Thus 
\begin{gather*}
\euler\left(\sum m_i G_i\right) > - 2^{6} t^2 \cdot \exp(14 |T_N^3|)\geq\\
\geq - 2^{6} t^2 \cdot \exp(14 \cdot 2^{9} t^2) > - \exp(2^{14} t^2).
\end{gather*}
The first inequality follows from Lemma~\ref{Lem:WeightOfFundamental}.
Recall that the Euler characteristic of $F$ is less than $-\exp(2^{15}
t^2)$. Hence some coefficient $n_i$ is non-zero, say $n_1 \geq 1$.
Rewrite the Haken sum in the form $F = J + K$ where $K = F_1$.  The
surface $J$ is then the sum of all the other terms.  Note that $\bdy K
= \emptyset$.

Properly isotope $J$ and $K$ to obtain the reduced sum $F = J' + K'$,
as defined in Section~\ref{Sec:Normal}.  Abuse notation again to
write $F = J + K$.

\section{The final annulus}
\label{Sec:LastAnnulus}

Here we analyze the patches, seams, and exchange bands of the reduced
sum $F = J + K$.  The exchange bands of this sum are our second source
of essential annuli and \Mobius~bands.  It is one of these which yields
the final contradiction.

\begin{claim}
\label{Clm:GoodSeamExists}
Every seam of $J + K$ is a closed regular neighborhood of an essential
simple closed curve in $F$.  Every exchange band is an annulus or
\Mobius~band.  Also, there is at least one seam which is not
boundary-parallel in $F$.
\end{claim}

\begin{proof}
The first two claims follows from Lemma~\ref{Lem:NoDiskPatches} and
from the fact that $\bdy K = \emptyset$.

To prove the third proceed by contradiction.  That is, suppose every
patch but at most one is an annulus.  As $\euler(K) < 0$ deduce that
$J$ is a union of annuli. Thus $\euler(J) = 0$.  It follows that
$\euler(F) = \euler(K) = \euler(F_1)$ which is impossible.
\end{proof}

\subsection{Examining the exchange band $A$}
Fix a seam of $J + K$ which is not boundary-parallel inside of $F$.
Let $A$ be the exchange band meeting this seam.  The argument finishes
with a discussion of this exchange band, $A$.  As above, the
hypothesis that $H$ is full rules out several possibilities.  The
eventual conclusion will be that the splitting $H$ admits a joined
annulus/annulus or \Mobius/annulus pair.  This implies that the
splitting has the DCP and yields the desired contradiction.

Let $\bdy_- A \subset \bdy A$ be the component of $\bdy A$ meeting the
seam which is not boundary-parallel in $F$. Let $B = \pi^{-1}(\bdy_-
A)$ be the vertical surface which is the union of all fibres meeting
$\bdy_- A$.  Recall, from Remark~\ref{Rem:Seams}, that all seams are
orientation-preserving in $F$.  Thus, $B$ is in fact an annulus.

\begin{claim}
\label{Clm:GoodSeamGivesEssentialAnnulus}
The annulus $B$ is essential in $V$. 
\end{claim}

\begin{proof}
Suppose that $B$ compresses in $V$. Then, by
Claim~\ref{Clm:HorizontalBdyIncompressibleInV}, the components of
$\bdy B$ bound disks in $\bdy_h Y$.  Thus $\bdy_- A = \pi(\bdy B)$ is
null-homotopic in $F$, contradicting Claim~\ref{Clm:GoodSeamExists}.

Suppose that $B$ is boundary-parallel in $V$. Let $C \subset H$ be
the annulus cobounded by $\bdy B$.  Let $X \subset V$ be the solid
torus bounded by $B \cup C$.  Note that $X$ admits a meridional disk
meeting each of $B$ and $C$ in a single arc.

If $C \subset \bdy_h Y$ then $X \subset Y$ and $X$ is the union of the
fibres of $Y$ meeting $C$, as $B$ is vertical. Thus $X$ is an
$I$--bundle over a \Mobius~band and any meridional disk meets $B$ in
at least two arcs, a contradiction.  So $C$ is not contained in
$\bdy_h Y$.  Thus $C$ contains some components of $\bdy (\bdy_h Y)$.  By
Claim~\ref{Clm:CornersOfY} every such is a core curve for $C$. It
follows that $\bdy B$ is boundary-parallel in $\bdy_h Y$ and thus
$\bdy_- A \subset B$ is boundary-parallel in $F$. This contradicts the
choice of exchange band $A$.
\end{proof}

\subsection{Constructing the annulus $A'$}

As $Y$ is an $I$--bundle there is a small isotopy of $Y$ to $Y'$ fixing
$F \cup \bdy N$ pointwise so that after the isotopy $A \cap Y' \subset
(\pi')^{-1}(\bdy A)$.  Here $\pi' \from Y' \to F$ is the obvious
bundle map.  We will also insist that $A \cap Y' = \neigh_A(\bdy A)$.
Again, abuse notation and refer to $Y'$ and $\pi'$ as $Y$ and $\pi$.


Let $A' = A \setminus \interior(Y)$.  Label $\bdy A'$ and $\bdy B$ so
that $\bdy_- A' = \bdy_+ B = A' \cap B$.  Isotope $H$, $\rel~Y \cup
R$, to minimize the quantity $|H \cap A'|$.  Once more, abuse notation
and denote the resulting Heegaard splitting surface by $H$.

\begin{claim}
\label{Clm:NoTrivialCurvesInA}
All curves of $A' \cap H$ are essential in $A'$.  
\end{claim}

\begin{proof}
Suppose $\delta \subset A' \cap H$ bounds an innermost disk $D \subset
A'$.  Now $\delta$ must be essential in $H$.  (If not $\delta$ bounds
a disk $E \subset H$.  As usual, because $\bdy_h Y$ is nonplanar, $E
\cap \bdy_h Y = \emptyset$.  So $E \cup D$ forms a two-sphere in $M$
and this bounds a ball, $U$ (Claim~\ref{Clm:MIrreducible}).  Since $(Y \cup
R) \cap U = \emptyset$ (Claim~\ref{Clm:YNotInBall}) there is an
isotopy of $H$, $\rel~Y \cup R$, which moves $E$ across $U$.  This
reduces $|A' \cap H|$, an impossibility.)  
Thus $D$ is an essential disk in $V$ or $W$ which is disjoint from
$\bdy_h Y$.  By Remark~\ref{Rem:AnnuliForY} the disk $D$ is one-half
of a disjoint disk/annulus pair, implying that $H$ has the DCP, a
contradiction.
\end{proof}


Let $A''$ be the closure of the component of $A' \setminus H$ such
that $A''$ contains $\bdy_- A' = \bdy_+ B$.  Note that $A'' \subset
W$.  Our task is to show that $A''$ and $B$ are a joined pair of
essential surfaces.  Begin as follows:

\begin{claim}
\label{Clm:AIsIncompressible}
If $A''$ is an annulus then $A''$ is incompressible in $W$.
\end{claim}

\begin{proof}
A compression of $A''$ would give a properly embedded disk $(D, \bdy D)
\subset (W, \bdy_+ B)$.  Recall that $B = \pi^{-1}(\bdy_- A)$, is
essential in $V$ (Claim~\ref{Clm:GoodSeamGivesEssentialAnnulus}).
Thus $D$ and $B$ form a joined disk/annulus pair which implies that
$H$ is weakly reducible.
\end{proof}

\begin{claim}
\label{Clm:FinalAnnulus}
The annulus or \Mobius~band $A''$ is essential in $W$.
\end{claim}

This will complete the proof of
Proposition~\ref{Prop:AlmostNormalPlusHighGenusImpliesDCP} as $A''$
and $B$ together form an joined annulus/annulus or \Mobius/annulus
pair.  This last shows that $H$ has the DCP and thus could not have
been a full Heegaard splitting.

\begin{proof}[Proof of Claim~\ref{Clm:FinalAnnulus}]
There is a penultimate dichotomy left to consider --- If $A''$ is
strictly contained in $A'$ then the situation is fairly simple.  If
$A'' = A'$ then a slightly more subtle analysis, using
Remark~\ref{Rem:ExchangeOnTunnel}, obtains.

Suppose that $\interior(A') \cap H \neq \emptyset$.  That is,
$A'' \neq A'$.  Then $A''$ is an annulus and hence
incompressible in $W$ by Claim~\ref{Clm:AIsIncompressible}.  Suppose that
$A''$ is boundary-parallel.  Thus $\bdy_+ B = \bdy_- A''$ is
parallel in $H$ to the curve $\bdy_+ A'' \subset (H \setminus
\bdy_h Y)$.  So $\bdy_+ B$ is boundary-parallel in the surface $\bdy_h
Y$ and thus $\bdy_- A \subset B$ is boundary-parallel in $F$.  This
contradicts the choice of seam made in Claim~\ref{Clm:GoodSeamExists}.

Suppose instead that $A' \cap H = \bdy A'$.  Then $A'' = A' \subset
W$.  There are two final cases: Either $A$ was an annulus or a
\Mobius~band.  Suppose that $A$ (and hence $A'$) is a \Mobius~band.
Now, $A'$ is properly embedded in the handlebody $W$ and $W$ is not
a solid torus.  Thus $A' = A''$ is essential, as desired. (See
Remark~\ref{Rem:Mobius}.)

At the last there is the possibility that $A$ is an annulus.  Recall
$Y$ has been isotoped slightly so that $A' = \closure{A \setminus Y}$
and $A' \subset W$.  Again, $A'' = A'$ is an annulus and hence
incompressible in $W$ by Claim~\ref{Clm:AIsIncompressible}.  Suppose
$A'$ is boundary-parallel in $W$.  Let $C \subset H$ be the annulus
cobounded by the two components of $\bdy A'$.  Let $X \subset W$ be
the solid torus bounded by $A' \cup C$.  If $C$ is not contained in
$\bdy_h Y$ then, as above, deduce that $\bdy_+ B$ is boundary-parallel
in $\bdy_h Y$.  This contradicts the choice of seam made in
Claim~\ref{Clm:GoodSeamExists}.  Conclude instead that $C \subset
\bdy_h Y$.

It follows $A$ is parallel, relative to $\bdy A$ and through $X$, into
the surface $F$.  As observed in Remark~\ref{Rem:ExchangeOnTunnel}
this gives a contradiction -- for the reader's convenience here are the
details: Let $C' \subset F$ be the annulus cobounded by the two
boundary components of $A$ and let $X'$ be the solid torus bounded by
$A \cup C'$.  Let $Z$ be the closure of $\neigh_{M \setminus F}(A)$.
Let $\double{A}$ be the closure of $\bdy Z \setminus F$.  Finally, let
$F'$ be a slight push-off of $(F \setminus Z) \cup \double{A}$.  In
the terminology of Jaco and Oertel~\cite{JacoOertel84} the surface
$F'$ is obtained by performing an {\em irregular exchange} along the
exchange annulus $A$.  Note that $F'$ is the union of a torus and a
surface, $F''$, isotopic to $F$.  For a schematic cross-sectional view
see Figure~\ref{Fig:ExchangeOnTunnel}.

\begin{figure}[ht!]\small
\psfrag{S}{$F$}
\psfrag{S''}{$F''$}
\psfrag{A}{$A$}  
\psfrag{X}{$X' \setminus Z$}
$$\begin{array}{c}
\epsfig{file=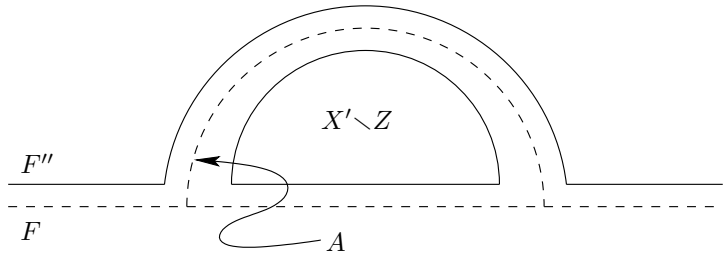, height = 3.3 cm}
\end{array}$$
\caption{An irregular exchange along $A$}
\label{Fig:ExchangeOnTunnel}
\end{figure}

Either $w(F'') < w(F)$ or $F''$ is not normal (consult
Figure~\ref{Fig:SmallNeighborhood}).  In the latter case there is a
weight reducing isotopy of $F''$ as in
Lemma~\ref{Lem:NormalizeIncomp}.  (See
Figure~\ref{Fig:BentArcWithIsotopy} with $F''$ replacing $F$.)  In
either case find a contradiction to the assumed lexicographic
minimality of $(w(\bdy F), w(F))$.  It follows that $A'$ was not
boundary-parallel in $W$.  Hence $A' = A''$ is essential in $W$.
\end{proof}

\section{Deducing the main theorem}
\label{Sec:MainTheorem}

Proposition~\ref{Prop:AlmostNormalPlusHighGenusImpliesDCP} is cast
into the realm of general splittings using
Theorem~\ref{Thm:NormalizeStrIrred}

\begin{theorem}
\label{Thm:HighGenusImpliesDCP}
Fix $M$, a closed, orientable three-manifold. There is a constant
$C(M)$ such that if $H \subset M$ is a Heegaard splitting, 
with $g(H) > C(M)$, then $H$ has the disjoint curve 
property.  
\end{theorem}

\begin{proof}
Let $T_M$ be a minimal triangulation of $M$ and take $C(M) =
\dcpconst$, where $t = |T_M^3|$.  Suppose that $g(H)$ is at least this
big.  If $H$ is strongly irreducible then, by
Theorem~\ref{Thm:NormalizeStrIrred}, isotope $H$ to be almost normal
with respect to $T_M$.  Thus $H$ has the DCP by
Proposition~\ref{Prop:AlmostNormalPlusHighGenusImpliesDCP}.  On
the other hand, if $H$ is weakly reducible then $H$ trivially
satisfies the DCP.
\end{proof}

\section{Conjecture}
\label{Sec:Conjecture}

Recall the generalized Waldhausen conjecture states that: 

\medskip{\bf Conjecture~\ref{Conj:GeneralizedWaldhausen}}\qua
{\sl 
If $M$ is closed, orientable, and atoroidal then $M$
contains only finitely many strongly irreducible splittings, up to
isotopy, in each genus.} 
\medskip

Now fix attention on a closed orientable manifold $M$.  If $M$ is a
lens space then all splittings are stabilizations of the genus one
splitting~\cite{BonahonOtal83}, and hence all splittings other than
the genus one splitting have the DCP.  We thus restrict ourselves to
manifolds with Heegaard genus two or larger.  If $M$ is toroidal then
by Thompson's Lemma~\ref{Lem:Thompson} all splittings in $M$ have the
DCP.  If $M$ is atoroidal then, by
Theorem~\ref{Thm:HighGenusImpliesDCP}, above a certain genus all
splittings have the DCP and, by
Conjecture~\ref{Conj:GeneralizedWaldhausen}~\cite{JacoRubinstein02},
there are only finitely many splittings up to isotopy below that
genus.  So it would follow that:

\begin{conjecture}
\label{Conj:OnlyFinManyFullSplittings}
In any closed, orientable three-manifold there are only fin\-itely many
full Heegaard splittings, up to isotopy.
\end{conjecture}

In addition to Conjecture~\ref{Conj:OnlyFinManyFullSplittings}
affirmative answers to any of the following questions would be very
interesting.

\begin{question}
Suppose $M$ is closed, orientable and $M$ admits a full Heegaard
splitting. By Lemmas~\ref{Lem:Haken}, \ref{Lem:Thompson},
and~\ref{Lem:Hempel} it follows that $M$ is irreducible, atoroidal,
and not a Seifert fibred space.  So from Thurston's geometrization
conjecture it would follow that $M$ is hyperbolic.  Weaker, but still
extremely interesting, would be a direct proof that $M$ had
word-hyperbolic, infinite, or simply non-trivial fundamental group.
\end{question}

Note that the converse statement, that every hyperbolic manifold
contains a full splitting, cannot be obtained.  It follows from work
of Rieck and Sedgwick~\cite{RieckSedgwick01}, combined with that of
Kobayashi~\cite{Kobayashi01}, that for ``most'' surgeries on
two-bridge knots the manifold obtained does {\em not} admit a full
splitting.  This answers a question of Luo's in the negative: There
are closed hyperbolic three-manifolds which do not admit a
``hyperbolic'' splitting. (See Section~5 of~\cite{Luo03}.)

More directly concerned with the work of this paper would be an
improvement of the bound given in Theorem~\ref{Thm:HighGenusImpliesDCP}.
A ``yes'' answer to the following question would be very satisfying:

\begin{question}
Do full Heegaard splittings always have minimal genus?
\end{question}

In light of the above, Sedgwick has asked a modified version of
Waldhausen's conjecture:

\begin{question}
If $M$ is non-Haken, does $M$ contain only finitely many
non-isotopic, irreducible Heegaard splittings?
\end{question}


Establishing this result would nicely complement
Conjecture~\ref{Conj:OnlyFinManyFullSplittings} as well as reaffirm the
Casson-Gordon program of studying non-Haken three-manifolds via their
strongly irreducible Heegaard splittings.

\newpage\appendix {\Large\bf Appendices}\small
\section{Surfaces}
\label{App:Surfaces}

Recall that the {\em complexity} of a compact surface $G$ is defined
to be the quantity $c(G) = \euler(G) + |\bdy G|$.  This is the same as
the Euler characteristic of the surface obtained by capping off all
boundary components of $G$ by disks. The lemma below essentially
states that Euler characteristic increases when a surface is
compressed.

\begin{lemma}
\label{Lem:SurfaceComplexityChange}
Suppose that $G$ and $F$ are compact, connected surfaces with $G
\subset F$. Then $c(F) \leq c(G)$.
\end{lemma}

\begin{proof}
If $G = F$ then the statement is trivial. Assume then that $\bdy G
\neq \emptyset$ and that $G$ has been isotoped so that $\bdy G \cap
\bdy F = \emptyset$.

Let $G' = \bigcup_{i = 1}^n G_i = F \setminus \interior(G)$, where
each $G_i$ is connected. It follows that $n \leq |\bdy G|$ and that
$|\bdy G'| = |\bdy G| + |\bdy F|$.  As $c(S) \leq 2$ for any compact,
connected surface $S$ we have:
$$ c(G') = \sum c(G_i) \leq 2n \leq 2 |\bdy G|.$$
Omitting the middle yields:
$$ \euler(G') + |\bdy G'| \leq 2 |\bdy G|.$$
Thus:
$$ \euler(G') + |\bdy F| \leq |\bdy G|.$$
Add $\euler(G)$ to both sides to obtain the desired inequality.
\end{proof}

\section{$I$--bundles}
\label{App:$I$-Bundles}

As a matter of terminology, a subset of an $I$--bundle is called {\em
vertical} if it is a union of fibres.  The next lemma supplies
us with many essential annuli inside of handlebodies.

\begin{lemma}
\label{Lem:NonSepVerticalAnnuliAreEssential}
Let $Y$ be an $I$--bundle embedded in a handlebody $V$ such that
$\bdy_h Y \subset \bdy V$.  Let $A$ be a properly embedded vertical
annulus of $Y$.  Suppose that $A$ is nonseparating inside of $Y$.
Then $A$ is an essential annulus in $V$.
\end{lemma}


\begin{proof}
As $A$ is nonseparating in $Y$, $A$ is nonseparating in $V$. Thus $A$
is not boundary-parallel. It is left to show that $A$ is incompressible.

Pick $B$, a vertical annulus or \Mobius~band in $Y$, such that $\alpha
= A \cap B$ is a single fibre.  For a contradiction suppose that $D
\subset V$ is a compressing disk for $A$. As $\bdy D$ is isotopic to
the core curve of $A$, $|\alpha \cap \bdy D|$ is odd.

However, by general position, $D \cap B$ is a compact one-manifold with
boundary $\alpha \cap \bdy D$. But compact one-manifolds have an even
number of points on their boundary. Thus $A$ is incompressible.
\end{proof}

\end{document}